\documentclass[12pt]{article}
\usepackage{graphicx}
\usepackage{amsmath}
\usepackage{graphicx}
\usepackage{amsfonts}
\usepackage{amssymb}
\usepackage{xcolor}
\usepackage{cite}
\usepackage{amsmath, amssymb ,amsthm, amsfonts, amsgen}
\DeclareGraphicsExtensions{.eps,.bmp,.jpg,.pdf,.mps,.png,.gif}
\numberwithin{equation}{section} 
\newcommand{\NN}{\mathbb{N}}

\newcommand{\RR}{\mathbb{R}}

\newtheorem{teo}{\sc Theorem}[section]
\newtheorem{prop}{\sc Proposition}[section]

\newtheorem{coro}{\sc Corollary}[section]

\newcommand{\dem}{\noindent {\it Proof.} \mbox{}}    
\newcommand{\lqqd}{\hfill $\square$}      

\begin{document}
\title{Junction of ferroelectric thin cylinders}

\author{
 Luisa Faella\footnote{Dipartimento di Ingegneria Elettrica e del l'Informazione ``M.Scarano",
Universit\`a degli Studi di Cassino e del Lazio Meridionale, Via G. Di Biasio n.43, 03043 Cassino (FR), Italy. E-mail: l.faella@unicas.it} \and Pedro Hern\'{a}ndez-Llanos\footnote{
 Instituto de Ciencias de la Ingenier\'ia, Universidad de O'Higgins, Avenida Libertador Bernardo O'Higgins 611, 2841935 Rancagua, Chile. E-mail: pedro.hernandez@uoh.cl {\bf(corresponding author)}}
 \and Ravi Prakash\footnote{Departamento de Matem\'aticas, Facultad de Ciencias F\'isicas y Matem\'aticas, Universidad de Concepci\'on, Avenida Esteban Iturra s/n, Barrio Universitario, Casilla 160-C, Concepci\'on, Chile. Email: rprakash@udec.cl}}
\date{\today}

\maketitle{}

\begin{abstract}

In this paper, starting from a non-convex and nonlocal $3D-$variational model for the electric polarization
in a ferroelectric material, and using an asymptotic process based on dimensional reduction, we
analyze junction phenomena for two vertical joined ferroelectric thin cilynders with differents small thicknesses. We obtain three
variational models for thin wire joined to a thin film, depending on the limit of the ratio between these two thicknesses.

{\it Key words and phrases.} Electric polarization, thin cylinders, junctions.

{\rm 2010 AMS Subject Classification. 35Q61 (primary); 78A25 (secondary).}
\end{abstract}

\section{Introduction}

There are materials with properties of great interest in the industry due to their various applications in the manufacture of electronic devices. Within the characteristics of such materials we can mention ferroelectricity, which is the property of to have a spontaneous electrical polarization that can be
reversed by the application of an external electric field, the behavior of
these materials is very similar to the one of ferromagnetic materials because hysteresis phenomena and Curie temperature appear $T_C$.  In these materials, for temperature above the Curie point, the internal energy has unique point named paraelectric phase (see \cite{Smith2005}).

\medskip

In the early 1920s, J. Valasek discovered the first ferroelectric material in Rochelle salt, a common salt but chemically and crystallographically complex. Later, another ferroelectric salt was discovered (KH$_2$PO$_4$) and in the decade of the 40s  ferroelectric properties were demonstrated in barium titanate (BaTiO$_3$) and lead titanate (PbTiO$_3$). For the history and applications of ferroelectric material, we refer to \cite{REF5, REF9, REF12} .

\medskip

Ferroelectric materials are widely used in different contexts, specifically in the manufacture of electronic circuits, integrated forms and storage devices. Among the most common, radio frequency identification cards (RFID) and the ferroelectric tunnel junction (FTJ) which appear in modern processors.

\medskip

In recent years some theoretical models for the study of ferroelectric phenomena were proposed. For instance, the mathematical modeling (in the static case) of thin structures of ferroelectric materials was studied starting from a non-convex and nonlocal $3D-$variational model for the electric polarization. Using asymptotic analysis based on dimensional reduction,
$2D-$variational models for thin films were obtained in \cite{REF22}, and $1D-$variational models for thin wires were obtained in \cite{REF23}. 

\medskip

 During the last decades the study of multistructures has aroused great interest in the mathematical community due to its applications in elasticity and high contrast materials, we cite to \cite{REF24, REF29, pedropaper2022, DemaioGaudielloSili2022}. For the limiting behaviour of $3D$ ferromagnetic problems in thin multistructures close to our model, we refer to \cite{AlCaLa, AmTo, Car, ChaFaPe, REF36bis, REF18,REF19, REF21,GoLuRoSl, sanchez,SlSo}. Recently, the fin junction of ferroelectric thin films and the junction of ferroelectric wires were examinated in \cite{carbonechaugaudie2018} and \cite{REF24*} respectively.

\medskip

\begin{figure}[h]
    \centering
    \includegraphics[width=0.7\textwidth]{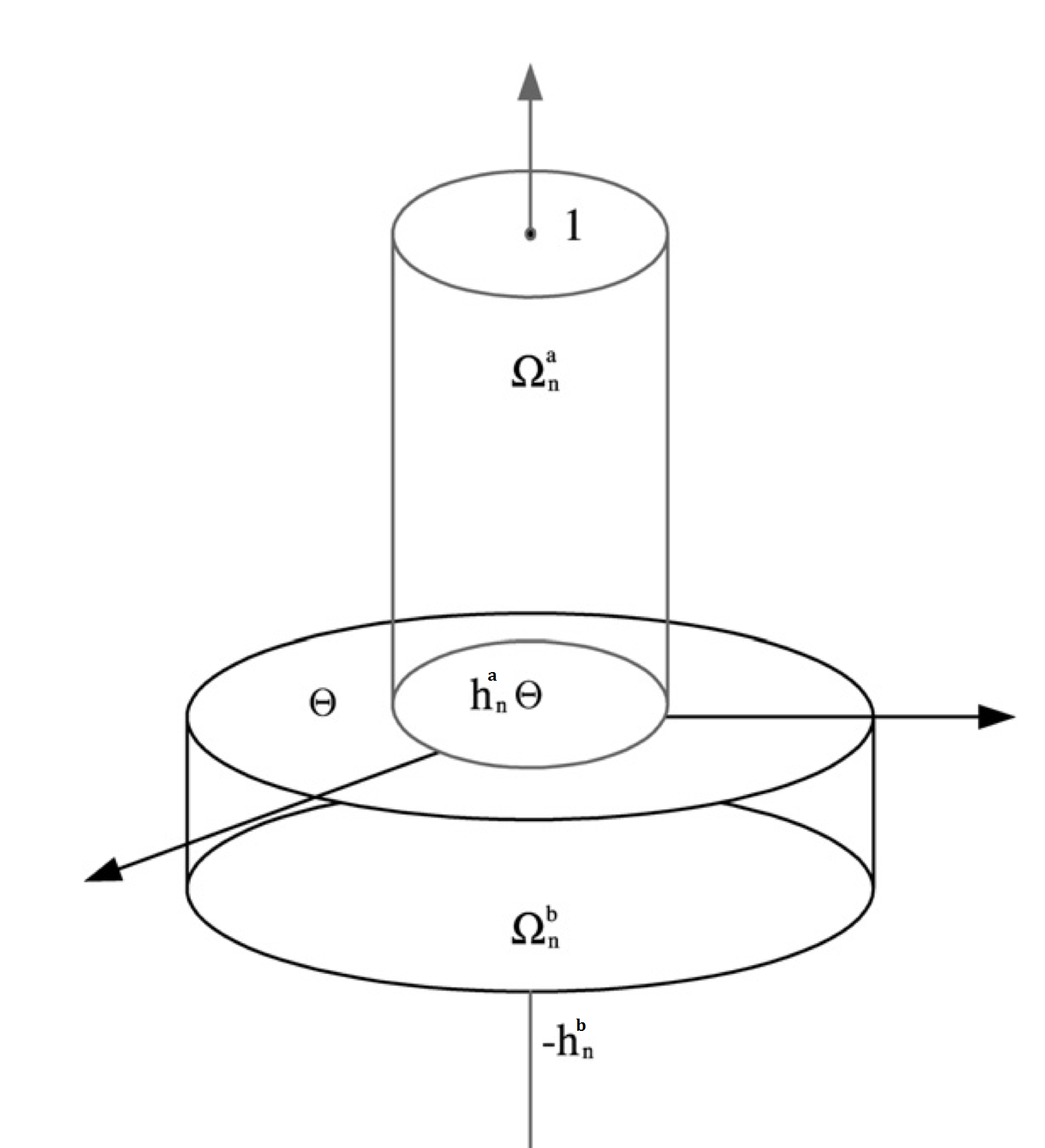}
    \caption{The set $\Omega_n$.}
    \label{fig:1}
\end{figure}

In this paper which was motivated by the structure given in \cite{REF19}, starting from the above described $3D-$variational model and using asymptotic method based on dimensional reduction, we analyze junction phenomena in a thin multistructure composed of two ferroelectric thin cylinders joined together (see Figure \ref{fig:1}). Such structure appears in some types of non-planar transistor used in the design of modern processors, the so-called Fin Field Effect Transistor (FinFET). The reduced model is justified for its simplicity and economical character from the computational point of view. We obtain three variational models for thin wire joined to a thin film, depending on how the reduction happens. In this study,
we do not consider any deformation of the ferroelectric material. 

\vspace{0.5 cm}

Let $\{h^a_{n}\}$ and $\{h^b_{n}\}\subset]0,1[$ be two sequences such that
\begin{equation}\label{eq:1}
\displaystyle\lim_{n}h^a_{n}=\lim_{n}h^b_{n}=0,\quad \lim_{n}\frac{h^b_{n}}{(h^a_{n})^2}=\ell\in[0,+\infty].
\end{equation} 

\medskip
For all $n\in\NN$, (see Figure \ref{fig:1}) set 
\begin{center}
$\Omega^{a}_{n}:=h^a_{n}\Theta\times[0,1[$, $\Omega^{b}_{n}:=\Theta\times\left]-h^b_{n}, 0\right[$,
\end{center}
and $\Omega_{n}:=\Omega^{a}_{n}\cup \Omega^{b}_{n}$, where $\Theta\subset\RR^2$ is  a bounded open connected set with smooth boundary such that $0':=(0,0)\in\Theta\subseteq ]-\frac{1}{2},\frac{1}{2}[\times]-\frac{1}{2},\frac{1}{2}[$.

\medskip
The multidomain $\Omega_{n}$  models a ferroelectric device which consist of two vertical cylinders $\Omega_{n}^{a}$ and $\Omega_{n}^{b}$ with one been placed upon the other. More precisely, the above cylinder have fixed height 1 and small cross section $h^a_n\Theta$. Whereas the belowan have small thickness $h^b_n$ and constant cross section $\Theta$. Here $h^a_n$ and $h^b_n$ are two small parameters converging to zero (see Figure \ref{fig:1}).

\vspace{0.5 cm}

We consider the following non-convex and non-local energy associated with $\Omega_{n}$ given in \cite{carbonechaugaudie2018,REF24*}:
\begin{equation}\label{eq:2}
\mathcal{E}_{n}:\textbf{P}\in\mathcal{P}_{n}\to \displaystyle\int_{\Omega_n}\left(\beta|\text{rot}\, \textbf{P}|^2+|\text{div}\, \textbf{P}|^2+\alpha\left(|\textbf{P}|^2-1\right)^2+|D\varphi_{\textbf{P}}|^2+(\textbf{F}_n\cdot\textbf{P})\right)dx,
\end{equation}
where,
\begin{equation}\label{rescaprobjuly19}
\mathcal{P}_n=\{\textbf{P}\in\left(H^1(\Omega_n)\right)^3:\textbf{P}\cdot\nu=0\quad\text{on}\,\partial\Omega_n\},
\end{equation}
$\textbf{P}$ is the spontaneous electric polarization field, $\alpha$ and $\beta$ are two positive constants, rot denotes the curl, $\textbf{F}_n\in\left(L^2(\Omega_n)\right)^3$, $\nu$ denotes the unit normal on $\partial\Omega_n$ and $\varphi_{\textbf{P}}\in H^1(\Omega_n)$ is the unique solution, up to an additive constant, of the boundary value problem
\begin{equation}\label{eq:3}
\begin{cases} 
      \text{div}\,\left(-\varepsilon_{0}D\varphi_{\textbf{P}}+\textbf{P}\right)=0 &\mbox{in }\Omega_n,   \\
      \left(-\varepsilon_0 D\varphi_{\textbf{P}}+\textbf{P}\right)\cdot\nu=0  &\mbox{on }\partial\Omega_n.
   \end{cases}
\end{equation}
Notice that $\textbf{F}_n$ is a normalization of external field. To see the essential characteristics of the model that we are considering, we can refer to \cite{REF5, REF9, REF12, REF30, REF31, REF32, REF33}.

 Since, 
\begin{equation}\label{eq:4}
||D\textbf{P}||^2_{\left(L^2(\Omega_n)\right)^9}=||\text{rot}\,\textbf{P}||^2_{\left(L^2(\Omega_n)\right)^3}+||\text{div}\,\textbf{P}||^2_{L^2(\Omega_n)},\quad\text{for all $\textbf{P}\in\mathcal{P}_n$ and all $n\in\NN$}
\end{equation}
(see Lemma 2.1 in \cite{REF23}), from the direct method of Calculus of Variations, the problem
\begin{equation}\label{eq:5}
\min \left\{\mathcal{E}_n(\textbf{P}):\textbf{P}\in\mathcal{P}_n\right\}
\end{equation}
admits a unique solution. Note that the problem (\ref{eq:5}) is an optimal control problem.
\medskip

The goal of this paper is to study the asymptotic behavior of $(\ref{eq:5})$ as $n$ diverges to infinity. We reformulate the problem using an appropriate rescaling (for more details see \cite{REF10}) in a fixed domain $\Omega=\Omega^a\cup \Omega^b$, where $\Omega^a:=\Theta\times]0,1[$ and $\Omega^{b}:=\Theta\times \left]-1, 0\right[$. Moreover, we rescale $\textbf{F}_n$ in $\Omega$ and we impose suitable convergence assumptions on these rescaled fields. Let  
\begin{align*}
E_0:q^a\in H^1(]0,1[)\mapsto \displaystyle &\nonumber |\Theta|\int_{0}^{1}\left(\left| \frac{d q^a}{d x_3}\right|^2+\alpha\left(|q^a|^2-1\right)^2+|D\psi^a_{(q^a, q^b)}|^2\right)dx_3,\\
&\nonumber\quad +\int_{0}^{1}\left(\int_{\Theta}f^a_3 dx'\cdot q^a\right)dx_3,\\
\end{align*}
\begin{align*}
E_{\infty}:q^b=(q_1^b,q_2^b)\in (H^1(\Theta))^2 \mapsto \displaystyle &\nonumber \int_{\Theta}\left(\beta\left|\text{rot}\, q^b\right|^2+\left|\text{div}\, q^b\right|^2+\alpha\left(|q^b|^2-1\right)^2+|D\psi^b_{(q^a, q^b)}|^2\right)dx',\\
&\nonumber\quad +\int_{\Theta}\left(\int_{-1}^{0}(f^b_1, f^b_2) dx_3\cdot q^b\right)dx',
\end{align*}

where
\begin{align*}
P=\left\{(q^a_3,q^b_1, q^b_2)\in H^1(]0, 1[)\times (H^1(\Theta))^2,\right.&\nonumber\, q^a_3(1)=0,\, (q^b_1,q^b_2,0)\cdot \nu^b=0\quad \text{on}\, \partial\Theta\\
&\left.\label{eq:16}q^a_3(0)=q^{b}_1(0')=q^{b}_2(0')=0 \right\}
\end{align*}
and $\left(\psi^{a}_{(q^a, q^b)}, \psi^{b}_{(q^a, q^b)}\right)$ is the unique solution of a suitable boundary value problem (see Section \ref{sect3}).
We study the asymptotic behavior of the rescaled energy $\mathcal{E}_n$ and we obtain a coupled $1D$ problem with $2D$ problem for $\ell\in(0,+\infty)$, a $1D$ problem for $\ell=0$ and a $2D$ problem for $\ell=+\infty$. More precisely, as the first problem with condition (\ref{eq:3}) we prove that:
\begin{itemize}
\item[(i)] For $\ell\in ]0,+\infty[$
\begin{equation}\label{july17limit1}\begin{array}{ll}
\displaystyle{\lim_n\min \left\{ \frac{1}{|\Omega_n|}{\mathcal E}_n(\textbf{P})\,:\,\textbf{P}\in\left(H^1(\Omega_n)\right)^3\,\,:\,\,\textbf{P}\cdot \nu=0\text{ on }\partial\Omega_n\right\}} \\\\
\displaystyle{=\frac{1}{2}\min\left\{E_0(q^a_3)\,:\,
q^a_3\in H_0^1(]0, 1[)\right\}}\\\\
\displaystyle{+\frac{\ell}{2}\min\left\{E_{\infty}(q^b_1,q^b_2):q^b_1, q^b_2\in  (H^1\left(\Theta\right))^2, (q^b_1,q^b_2,0)\cdot \nu^b=0\,\, \text{on}\, \partial\Theta, q^{b}_1(0')=q^{b}_2(0')=0 \right\}},
\end{array}
\end{equation}
more precisely, Theorem \ref{teo:1}. Where the junction condition is proved 
\begin{equation*}
q^a_3(0)=q^{b}_1(0')=q^{b}_2(0')
\end{equation*}
and play a important role in our results.
\item[(ii)] For $\ell=0$
\begin{equation}\label{july17limit2}\begin{array}{ll}
\displaystyle{\lim_n\min \left\{ \frac{1}{|\Omega_n|}{\mathcal E}_n(\textbf{P})\,:\,\textbf{P}\in\left(H^1(\Omega_n)\right)^3\,\,:\,\,\textbf{P}\cdot \nu=0\text{ on }\partial\Omega_n\right\}} \\\\
\displaystyle{=\frac{1}{2}\min\left\{E_0(q^a_3)\,:\,
q^a_3\in H_0^1(]0, 1[)\right\}},
\end{array}
\end{equation}
more precisely, Theorem \ref{teo:2}.
\item[$(iii)$] For $\ell=+\infty$
\begin{equation}\label{july17limit3}\begin{array}{ll}
\displaystyle{\lim_n\min \left\{ \frac{1}{|\Omega_n|}{\mathcal E}_n(\textbf{P})\,:\,\textbf{P}\in\left(H^1(\Omega_n)\right)^3\,\,:\,\,\textbf{P}\cdot \nu=0\text{ on }\partial\Omega_n\right\}} \\\\
\displaystyle{=\frac{1}{2}\min\left\{E_{\infty}(q^b_1,q^b_2): (q^b_1, q^b_2)\in  (H^1\left(\Theta\right))^2, (q^b_1,q^b_2,0)\cdot \nu^b=0\,\, \text{on}\, \partial\Theta, q^{b}_1(0')=q^{b}_2(0')=0 \right\}},
\end{array}
\end{equation}
more precisely, Theorem \ref{teo:3}.
\end{itemize}
Also, we show analogous results if we consider an energy like the following
\begin{equation*}
\mathcal{S}_{n}:\textbf{P}\in\left(H^1(\Omega_n)\right)^3\to \displaystyle\int_{\Omega_n}\left(|D\textbf{P}|^2+\alpha\left(|\textbf{P}|^2-1\right)^2+|D\varphi_{\textbf{P}}|^2+(\textbf{F}_n\cdot\textbf{P})\right)dx,
\end{equation*}
with the problems
\begin{equation*}
\min \left\{\mathcal{S}_n(\textbf{P}):\textbf{P}// e_3\,\text{on}\,\partial\Omega_n\right\}, 
\end{equation*}
where $\textbf{P}$ is the electric polarization, $e_3=(0,0,1)$ and $//$ is the symbol of parallelism. For details, see Section \ref{chap:6}.

\vspace{0.5 cm}

The paper is organized as follows: In Section \ref{sect2} we give the notation and the rescaled version for problem (\ref{eq:5}). In Section \ref{sect3}, we show the main results Theorem \ref{teo:1}, Theorem \ref{teo:2} and Theorem \ref{teo:3} which corresponds to our limit models for for the cases $\ell\in]0,\infty[$, $\ell=0$ and $\ell=+\infty$ respectively. Section \ref{cha:4} and Section \ref{cha:5} are devoted to prove the Theorem \ref{teo:1} and Theorem \ref{teo:2} respectively. Finally, in Section \ref{chap:6} we introduce a important choice for energetic approach  $\mathcal{S}_n$ which consider the energetic contribution for the polarization field given by the gradient term. We show some results for $\mathcal{S}_{n}$ with the boundary condition $\textbf{P}// e_3$ and $\textbf{P}\cdot \nu=0$ on $\partial\Omega_n$ for all the cases of $\ell$ and the relationship between the limiting energies for $\mathcal{S}_n$ and $\mathcal{E}_n$.

\vspace{1.0 cm}

\section{The rescaled problem}\label{sect2}

As it is usual (see \cite{REF10}) problem (\ref{eq:5}) is reformulated on a fixed domain through the maps
\begin{equation}\label{eq:6}
\begin{cases} 
      x=(x_1,x_2,x_3)\in \Omega^a=\Theta\times]0,1[\mapsto (h^a_nx', x_3)\in\text{Int}(\Omega^a_n), &   \\
      x=(x_1,x_2,x_3)\in \Omega^b=\Theta\times \left]-1, 0\right[\mapsto (x', h^b_nx_3)\in\Omega^b_n, &
   \end{cases}
\end{equation}
where $x'=(x_1,x_2)$ and $\text{Int}(\Omega^a_n)$ denotes the interior of $\Omega^a_n$. Precisely, for every $n\in\NN$ set
\begin{align*}
D^a_n: p^a\in\left(H^1(\Omega^a)\right)^k &\mapsto \displaystyle\left(\frac{1}{h^a_n}\frac{\partial p^a}{\partial x_1}, \frac{1}{h^a_n}\frac{\partial p^a}{\partial x_2}, \frac{\partial p^a}{\partial x_3}\right)\in \left(L^2(\Omega^a)\right)^{3k},\quad k\in\{1,3\},\\
D^b_n: p^b\in\left(H^1(\Omega^b)\right)^k &\mapsto \displaystyle\left(\frac{\partial p^b}{\partial x_1}, \frac{\partial p^b}{\partial x_2}, \frac{1}{h^b_n}\frac{\partial p^b}{\partial x_3}\right)\in \left(L^2(\Omega^b)\right)^{3k},\quad k\in\{1,3\},
\end{align*}

\begin{align*}
\text{div}^a_n: p^a=(p^a_1, p^a_2, p^a_3)\in\left(H^1(\Omega^a)\right)^3 &\mapsto \displaystyle\frac{1}{h^a_n}\frac{\partial p^a_1}{\partial x_1}+ \frac{1}{h^a_n}\frac{\partial p^a_2}{\partial x_2}+\frac{\partial p^a_3}{\partial x_3} \in L^2(\Omega^a),\\
\text{div}^b_n: p^b=(p^b_1, p^b_2, p^b_3)\in\left(H^1(\Omega^b)\right)^3 &\mapsto \displaystyle\frac{\partial p^b_1}{\partial x_1}+\frac{\partial p^b_2}{\partial x_2}+\frac{1}{h^b_n}\frac{\partial p^b_3}{\partial x_3} \in L^2(\Omega^b),
\end{align*}
\begin{align*}
\text{rot}^a_n: p^a=(p^a_1, p^a_2, p^a_3)\in\left(H^1(\Omega^a)\right)^3 &\mapsto \displaystyle\left(\frac{1}{h^a_n}\frac{\partial p^a_3}{\partial x_2}-\frac{\partial p^a_2}{\partial x_3}, \frac{\partial p^a_1}{\partial x_3}-\frac{1}{h^a_n}\frac{\partial p^a_3}{\partial x_1}, \frac{1}{h^a_n}\frac{\partial p^a_2}{\partial x_1}-\frac{1}{h^a_n}\frac{\partial p^a_1}{\partial x_2}\right)\in \left(L^2(\Omega^a)\right)^{3},\\
\text{rot}^b_n: p^b=(p^b_1, p^b_2, p^b_3)\in\left(H^1(\Omega^a)\right)^3 &\mapsto  \displaystyle\left(\frac{\partial p^b_3}{\partial x_2}-\frac{1}{h^b_n}\frac{\partial p^b_2}{\partial x_3}, \frac{1}{h^b_n}\frac{\partial p^b_1}{\partial x_3}-\frac{\partial p^b_3}{\partial x_1}, \frac{\partial p^b_2}{\partial x_1}-\frac{\partial p^b_1}{\partial x_2}\right)\in \left(L^2(\Omega^b)\right)^{3},
\end{align*}
and
\begin{equation}\label{eq:7}
\begin{cases} 
      f^a_n: x=(x_1,x_2,x_3)\in \Omega^a\mapsto \textbf{F}_n(h_n^ax', x_3), &   \\
      f^b_n: x=(x_1,x_2,x_3)\in \Omega^b\mapsto \textbf{F}_n(x', h^b_nx_3), &
   \end{cases}
\end{equation}
and define the set
\begin{align}
P_n=\left\{(p^a, p^b)\in\left(H^1(\Omega^a)\right)^3\times \left(H^1(\Omega^b\right))^3\right. &\nonumber: p^a\cdot \nu^a=0\quad \text{on}\, \partial\Omega^a \setminus(\Theta\times\{0\}) \\
&\nonumber \left. p^b\cdot \nu^b=0\quad \text{on}\, \partial\Omega^b \setminus(\Theta\times\{0\})\right.\\
&\nonumber \left. p^b_3=0\quad \text{on}\, \left(\Theta\setminus h^a_n\Theta\right)\times\{0\}\right.\\
&\label{eq:8} \left. p^a(x', 0)=p^b(h^a_n x', 0)\quad \text{a.e. in}\, \Theta \right\},
\end{align}
where $\nu^a$ and $\nu^b$ denote de unit outer normals on $\partial\Omega^a$ and $\partial\Omega^b$, respectively, and 
\begin{equation}\label{eq:9}
U_n=\{(\phi^a, \phi^b)\in H^1(\Omega^a)\times H^1(\Omega^b): \phi^a (x', 0)=\phi^b(h^a_n x',0)\quad \text{ for $x'$ a.e. in}\, \Theta\}.
\end{equation}
Then $\mathcal{E}_n$ defined in (\ref{eq:2}) is rescaled by
\begin{align}
E_{n}:(p^a, p^b)\in P_{n}\to \displaystyle &\nonumber (h^a_n)^2\int_{\Omega^a}\left(\beta|\text{rot}_n^a\, p^a|^2+|\text{div}_n^a\, p^a|^2+\alpha\left(|p^a|^2-1\right)^2-|D_n^a\phi^a_{(p^a, p^b)}|^2+(f^a_n\cdot p^a)\right)dx\\
&\label{eq:10}\quad +h^b_n\int_{\Omega^b}\left(\beta|\text{rot}_n^b\, p^b|^2+|\text{div}_n^b\, p^b|^2+\alpha\left(|p^b|^2-1\right)^2-|D_n^b\phi^b_{(p^a, p^b)}|^2+(f^b_n\cdot p^b)\right)dx,
\end{align}
where $\left(\phi^a_{(p^a, p^b)},\phi^b_{(p^a, p^b)} \right)$ is the unique solution of
\begin{equation}\label{eq:11}
\begin{cases} 
      \left(\phi^a_{(p^a, p^b)},\phi^b_{(p^a, p^b)} \right)\in U_n, \quad \displaystyle\int_{\Omega^a}\phi^a_{(p^a, p^b)}dx=0, &   \\
      (h^a_n)^2\displaystyle\int_{\Omega^a}\left(\left(-D^a_n\phi^a_{(p^a, p^b)}+p^a\right)\cdot D^a_n \phi^a\right)dx+h^b_n\displaystyle\int_{\Omega^b}\left(\left(-D^b_n\phi^b_{(p^a, p^b)}+p^b\right)\cdot D^b_n \phi^b\right)dx=0, & 
   \end{cases}
\end{equation}
 $(\phi^a, \phi^b)\in U_n$,
which rescales a weak formulation of (\ref{eq:3}), \textit{i.e.}
\begin{equation}\label{eq:12}
\varphi_{\textbf{P}}\in H^1(\Omega_n),\quad \displaystyle\int_{\Omega^a_n}\varphi_{\textbf{P}}dx=0,\quad \displaystyle\int_{\Omega_n}\left(\left(-\varepsilon_{0}D\varphi_{\textbf{P}}+\textbf{P}\right)\cdot D\varphi\right)dx=0,\quad \varphi\in H^1(\Omega_n).
\end{equation}
The Lax-Milgram Theorem provides that (\ref{eq:12}) admits solutions and it is unique. Note that if $\textbf{P}_n$ solves (\ref{eq:5}), then $(p^a_n, p^b_n)$ defined by
\begin{equation*}
p^a_n(x_1, x_2, x_3)=\textbf{P}_n(h^a_n x', x_3)\quad \text{in}\, \Omega^a,\quad p^b_n(x_1, x_2, x_3)=\textbf{P}_n(x', h^b_n x_3)\quad \text{in}\, \Omega^b,
\end{equation*}
solves
\begin{equation}\label{eq:13}
\min \{E_n((p^a, p^b)): (p^a, p^b)\in P_n\}.
\end{equation}
Conversely, if $(p^a_n, p^b_n)$ solves (\ref{eq:13}), then $\textbf{P}_n$ defined by
\begin{equation*}
\textbf{P}_n(x_1, x_2, x_3)=p^a_n\left(\frac{x'}{h^a_n}, x_3\right)\, \text{in}\quad \Omega^a_n,\quad \textbf{P}_n(x_1, x_2, x_3)=p^a_n\left(x', \frac{x_3}{h^b_n}\right)\, \text{in}\quad \Omega^b_n,
\end{equation*}
solves $(\ref{eq:5})$. Therefore, the goal of this paper becomes to study the asymptotic behavior, as $n$ diverges of (\ref{eq:13}). To this end, we assume
\begin{equation}\label{eq:14}
\begin{array}{rcl}
          f^a_n\rightharpoonup f^a & \text{weakly in} & \left(L^2(\Omega^a)\right)^3, \\
          f^b_n\rightharpoonup f^b & \text{weakly in} & \left(L^2(\Omega^b)\right)^3.
\end{array}
\end{equation}

\medskip
Note that rescalings in (\ref{eq:6}) transform (\ref{eq:4}) into
\begin{align}
(h^a_n)^2||D^a_n p^a||^2_{\left(L^2(\Omega^a)\right)^9}+h^b_n||D^b_n p^b||^2_{\left(L^2(\Omega^b)\right)^9}&\nonumber=(h^a_n)^2||\text{rot}\,_n^a p^a||^2_{\left(L^2(\Omega^a)\right)^3}+(h^a_n)^2||\text{div}\,^a_n p^a||^2_{L^2(\Omega^a)}\\
&\label{eq:15}\quad+h^b_n||\text{rot}\,_n^b p^b||^2_{\left(L^2(\Omega^b)\right)^3}+h^b_n||\text{div}\,^b_n p^b||^2_{L^2(\Omega^b)},
\end{align}
for all $(p^a, p^b)\in P_n$ and all $n\in\NN$.

\section{The main results}\label{sect3}

\subsection{The case $\ell\in]0,+\infty[$}
\medskip
This subsection is devoted to starting Theorem \ref{teo:1} describing the limit behavior of (\ref{eq:13}) when (\ref{eq:1}) is assumed with $\ell\in]0, +\infty[$. Theorem \ref{teo:1} will be proved in Section \ref{cha:4}.

\vspace{0.5 cm}
Set
\begin{align}
P=\left\{(q^a_3,q^b_1, q^b_2)\in H^1(]0, 1[)\times (H^1(\Theta))^2,\right.&\nonumber\, q^a_3(1)=0,\, (q^b_1,q^b_2,0)\cdot \nu^b=0\quad \text{on}\, \partial\Theta\\
&\left.\label{eq:16}q^a_3(0)=q^{b}_1(0')=q^{b}_2(0')=0 \right\},
\end{align}
\begin{equation}\label{eq:17}
U=\left\{(\psi^a, \psi^b)\in \left(H^1(]0, 1[)\right)\times (H^1(\Theta)),\quad\psi^a(0)=\psi^{b}(0') \right\}.
\end{equation}
and
\begin{align}
E:(q^a, q^b)\in P\mapsto \displaystyle &\nonumber |\Theta|\int_{0}^{1}\left(\left| \frac{d q^a}{d x_3}\right|^2+\alpha\left(|q^a|^2-1\right)^2+|D\psi^a_{(q^a, q^b)}|^2\right)dx_3,\\
&\nonumber\quad +\int_{0}^{1}\left(\int_{\Theta}f^a_3 dx'\cdot q^a\right)dx_3,\\
&\nonumber\quad +\ell\int_{\Theta}\left(\beta\left|\text{rot}\, q^b\right|^2+\left|\text{div}\, q^b\right|^2+\alpha\left(|q^b|^2-1\right)^2+|D\psi^b_{(q^a, q^b)}|^2\right)dx',\\
&\label{eq:18}\quad +\ell\int_{\Theta}\left(\int_{-1}^{0}(f^b_1, f^b_2) dx_3\cdot q^b\right)dx',
\end{align}
where $f^a_3$ is the third component of $f^a$ and $f^b_1$ and  $f^b_2$ are the first and second component of $f^b$ with $f^a, f^b$ are defined in (\ref{eq:14}), and $\left(\psi^{a}_{(q^a, q^b)}, \psi^{b}_{(q^a, q^b)}\right)$ is the unique solution of 
\begin{equation}\label{eq:19}
\begin{cases} 
      \left(\psi^a_{(q^a, q^b)},\psi^b_{(p^a, p^b)} \right)\in U, \quad \displaystyle\int_{0}^{1}\psi^a_{(q^a, q^b)}dx=0, &   \\
      |\Theta|\displaystyle\int_{0}^{1}\left(\left(-D_{x_3}\psi^a_{(q^a, q^b)}+q^a\right)\cdot D_{x_3} \psi^a\right)dx_3+\ell\displaystyle\int_{\Theta}\left(\left(-D_{x'}\psi^b_{(q^a, q^b)}+q^b\right)\cdot D_{x'} \psi^b\right)dx'=0, & 
 \end{cases}
\end{equation}
for $(\psi^a, \psi^b)\in U$.
Note that $(\ref{eq:19})$ admits a solution and it is unique since the set
\begin{equation*}
\left\{(\psi^a, \psi^b)\in U:\quad \displaystyle\int_{0}^{1}\psi^a dx_3=0\right\},
\end{equation*}
is a Hilbert space with the inner product
\begin{equation*}
\left<(\psi^a, \psi^b), (\varphi^a, \varphi^b)\right>=\displaystyle|\Theta|\int_{0}^{1}D\psi^a D\varphi^a dx_3+\ell\int_{\Theta}D\psi^b D\varphi^b dx'.
\end{equation*}

\vspace{1.0 cm}

\begin{teo}\label{teo:1}
Assume (\ref{eq:1}) with $\ell\in]0,+\infty[$ and (\ref{eq:14}). For every $n\in \NN$, let $(p^a_n, p^b_n)$ a solution of $(\ref{eq:13})$, and $\left(\phi^a_{(p^a_n, p^b_n)}, \phi^b_{(p^a_n, p^b_n)}\right)$ be the unique solution of $(\ref{eq:11})$ with $(p^a, p^b)=(p^a_n, p^b_n)$. Moreover, let $P$ and $E$ defined by (\ref{eq:16}) and (\ref{eq:18}), (\ref{eq:19}) respectively. Then there exist an increasing sequence of positive integer numbers $\{n_i\}_{i\in\NN}$ and (in possible dependence on the subsequence) $(\hat{p}^a, \hat{p}^b_1, \hat{p}^b_2)\in P$ such that
\begin{equation}\label{eq:20}
\begin{array}{rcl}
          p^a_{n_i}\rightarrow (0, 0, \hat{p}^a) & \text{strongly in} & \left(H^1(\Omega^a)\right)^3\,  \text{ and strongly in}\,  \left(L^4(\Omega^a)\right)^3,\\
          p^b_{n_i}\rightarrow (\hat{p}^b_1, \hat{p}^b_2 , 0) & \text{strongly in} & \left(H^1(\Omega^b)\right)^3\,  \text{ and strongly in}\,  \left(L^4(\Omega^b)\right)^3,
\end{array}
\end{equation}

\begin{equation}\label{eq:21}
\begin{cases} 
      \displaystyle\left(\frac{1}{h^a_n}\frac{\partial p^a_n}{\partial x_1}, \frac{1}{h^a_n}\frac{\partial p^a_n}{\partial x_2}\right)\rightarrow (0,0) & \text{strongly in}\, \left(L^2(\Omega^a)\right)^{3}\times \left(L^2(\Omega^a)\right)^{3}, \\
            \displaystyle \frac{1}{h^b_n}\frac{\partial p^b_n}{\partial x_3}\rightarrow 0 & \text{strongly in}\, \left(L^2(\Omega^b)\right)^{3}, 
 \end{cases}
\end{equation}

\begin{equation}\label{eq:22}
\begin{cases} 
      \displaystyle\left(\phi^a_{(p^a_{n_i},p^b_{n_i})}, \phi^b_{(p^a_{n_i},p^b_{n_i})} \right)\rightarrow(\psi^a_{(\hat{p}^a, \hat{p}^b)}, \psi^b_{(\hat{p}^a, \hat{p}^b)}) & \text{strongly in}\, H^1(\Omega^a)\times H^1(\Omega^b), \\
     \displaystyle\left(\frac{1}{h^a_n}\frac{\partial \phi^a_{(p^a_{n_i},p^b_{n_i})}}{\partial x_1} ,\frac{1}{h^a_n}\frac{\partial \phi^a_{(p^a_{n_i},p^b_{n_i})}}{\partial x_2} \right)\rightarrow(0,0) & \text{strongly in}\, \left(L^2(\Omega^a)\right)^2,\\
       \displaystyle\frac{1}{h^b_n}\frac{\partial  \phi^b_{(p^a_{n_i},p^b_{n_i})} }{\partial x_3}\rightarrow 0 & \text{strongly in}\, \left(L^2(\Omega^b)\right),  
 \end{cases}
\end{equation}
where $(\hat{p}^a, (\hat{p}^b_1, \hat{p}^b_2))$ solves
\begin{equation}\label{eq:23}
E((0,0,\hat{p}^a), (\hat{p}^b_1, \hat{p}^b_2,0))=\min \{E((0,0,q^a_3), (q^b_1, q^b_2,0)): (q^a_3, q^b_1,q^b_2)\in P\},
\end{equation}
 and $\left(\psi^a_{(\hat{p}^a, \hat{p}^b)}, \psi^b_{(\hat{p}^a, \hat{p}^b)}\right)$ is the unique solution of (\ref{eq:19}). Moreover
 \begin{equation}\label{eq:24}
 \displaystyle\lim_n \frac{E_n((p^a_n, p^b_n))}{(h^a_n)^2}=E((0,0,\hat{p}^a), (\hat{p}^b_1, \hat{p}^b_2,0)).
\end{equation}
\end{teo}

\vspace{0.5 cm}
\subsection{The case $\ell=0$}
This subsection is devoted to starting Theorem \ref{teo:2} describing the limit behavior of (\ref{eq:13}) when (\ref{eq:1}) is assumed with $\ell=0$. Theorem \ref{teo:2} will be proved in Section 5.

\vspace{0.5 cm}
Set
\begin{equation}\label{eq:16*new}
P_0=\left\{q^a\in H^1(]0, 1[),q^a(0)=q^{a}(1)=0 \right\}=H^1_0(]0,1[)
\end{equation}
and
\begin{align}
E_0: q^a\in H^1(]0,1[)\mapsto \displaystyle &\nonumber |\Theta|\int_{0}^{1}\left(\left| \frac{d q^a}{d x_3}\right|^2+\alpha\left(|q^a|^2-1\right)^2+|D\psi^a_{(q^a, q^b)}|^2\right)dx_3,\\
&\label{eq:18*}\quad +\int_{0}^{1}\left(\int_{\Theta}f^a_3 dx_1dx_2\cdot q^a\right)dx_3,
\end{align}
where $f^a_3$ and $f^b_1$ are defined in (\ref{eq:14}), and $\psi^{a}_{q^a}$ is the unique solution of 
\begin{equation}\label{eq:19*}
\begin{cases} 
      \psi^a_{q^a}\in U, \quad \displaystyle\int_{0}^{1}\psi^a_{q^a}dx=0, &   \\
      |\Theta|\displaystyle\int_{0}^{1}\left(\left(-D_{x_3}\psi^a_{q^a}+q^a\right)\cdot D_{x_3} \psi^a\right)dx_3=0, & 
 \end{cases}
\end{equation}
for $\psi^a\in H^1(]0,1[)$.

\vspace{1.0 cm}

\begin{teo}\label{teo:2}
Assume (\ref{eq:1}) with $\ell=0$ and (\ref{eq:14}). For every $n\in \NN$, let $(p^a_n, p^b_n)$ a solution of $(\ref{eq:13})$, and $\left(\phi^a_{(p^a_n, p^b_n)}, \phi^b_{(p^a_n, p^b_n)}\right)$ be the unique solution of $(\ref{eq:11})$ with $(p^a, p^b)=(p^a_n, p^b_n)$. Moreover, let $P_0$ and $E_0$ defined by (\ref{eq:16*new}) and (\ref{eq:18*}), (\ref{eq:19*}) respectively. Then there exist an increasing sequence of positive integer numbers $\{n_i\}_{i\in\NN}$ and (in possible dependence on the subsequence) $\hat{p}^a\in P_0$ such that
\begin{equation}\label{eq:20*}
\begin{array}{rcl}
          p^a_{n_i}\rightarrow (0, 0, \hat{p}^a) & \text{strongly in} & \left(H^1(\Omega^a)\right)^3\,  \text{ and strongly in}\,  \left(L^4(\Omega^a)\right)^3,\\
          \displaystyle\frac{\sqrt{h^b_n}}{h^a_n}p^b_{n_i}\rightarrow 0 & \text{strongly in} & \left(H^1(\Omega^b)\right)^3\,  \text{ and strongly in}\,  \left(L^4(\Omega^b)\right)^3,
\end{array}
\end{equation}

\begin{equation}\label{eq:21*}
\begin{cases} 
      \displaystyle\left(\frac{1}{h^a_n}\frac{\partial p^a_n}{\partial x_1}, \frac{1}{h^a_n}\frac{\partial p^a_n}{\partial x_2}\right)\rightarrow(0,0) & \text{strongly in}\, \left(L^2(\Omega^a)\right)^{3}\times \left(L^2(\Omega^a)\right)^{3}, \\
      \displaystyle\frac{1}{\sqrt{h^b_n}h^a_n}\frac{\partial p^b_n}{\partial x_3} \rightarrow 0 & \text{strongly in}\, \left(L^2(\Omega^b)\right)^{3}, 
 \end{cases}
\end{equation}

\begin{equation}\label{eq:22*}
\begin{cases} 
      \displaystyle\left(\phi^a_{(p^a_{n_i},p^b_{n_i})}, \frac{\sqrt{h^b_n}}{h^a_n}\phi^b_{(p^a_{n_i},p^b_{n_i})} \right)\rightarrow(\psi^a_{\hat{p}^a}, 0) & \text{strongly in}\, \left(H^1(\Omega^a)\right)\times \left(H^1(\Omega^b)\right), \\
     \displaystyle\left(\frac{1}{h^a_n}\frac{\partial \phi^a_{(p^a_{n_i},p^b_{n_i})}}{\partial x_1} ,\frac{1}{h^a_n}\frac{\partial \phi^a_{(p^a_{n_i},p^b_{n_i})}}{\partial x_2} \right)\rightarrow(0,0) & \text{strongly in}\, \left(L^2(\Omega^a)\right)^2,\\
       \displaystyle\frac{1}{\sqrt{h^b_n}h^a_n}\frac{\partial  \phi^b_{(p^a_{n_i},p^b_{n_i})} }{\partial x_3}\rightarrow 0& \text{strongly in}\, \left(L^2(\Omega^b)\right),  
 \end{cases}
\end{equation}
where $\hat{p}^a$ solves
\begin{equation}\label{eq:23*}
E_0(\hat{p}^a)=\min \{E_0(q^a): q^a\in P_0\},
\end{equation}
 and $\psi^a_{\hat{p}^a}$ is the unique solution of (\ref{eq:19*}) with $q^a=\hat{p}^a$. Moreover
 \begin{equation}\label{eq:24*}
 \displaystyle\lim_n \frac{E_n((p^a_n, p^b_n))}{(h^a_n)^2}=E_0(\hat{p}^a).
\end{equation}
\end{teo}

\vspace{0.5 cm}

\subsection{The case $\ell=+\infty$}
\medskip
This subsection is devoted to starting Theorem \ref{teo:3} describing the limit behavior of (\ref{eq:13}) when (\ref{eq:1}) is assumed with $\ell=+\infty$ and $h^b_n\ll \sqrt{h^a_n}$. Here we assume that the function $\left(\phi^a_{(p^a, p^b)},\phi^b_{(p^a, p^b)} \right)$ involved in (\ref{eq:11}) is the unique solution of the following problem:
\begin{equation}\label{eq:11*}
\begin{cases} 
      \left(\phi^a_{(p^a, p^b)},\phi^b_{(p^a, p^b)} \right)\in U_n, \quad \displaystyle\int_{\Omega^b}\phi^b_{(p^a, p^b)}dx=0, &   \\
      (h^a_n)^2\displaystyle\int_{\Omega^a}\left(\left(-D^a_n\phi^a_{(p^a, p^b)}+p^a\right)\cdot D^a_n \phi^a\right)dx+(h^b_n)^2\displaystyle\int_{\Omega^b}\left(\left(-D^b_n\phi^b_{(p^a, p^b)}+p^b\right)\cdot D^b_n \phi^b\right)dx=0, & 
   \end{cases}
\end{equation}
 $(\phi^a, \phi^b)\in U_n$,
i.e., the assumption $\int_{\Omega^a}\phi^a_{(p^a, p^b)}dx=0$ is replaced with 
\begin{equation*}
\int_{\Omega^b}\phi^b_{(p^a, p^b)}dx=0, 
\end{equation*}
or equivalently, in (\ref{eq:12}) the assumption $\int_{\Omega^a_n}\varphi_{\textbf{P}}dx=0$ is replaced with
\begin{equation*}
\int_{\Omega^b_n}\varphi_{\textbf{P}}dx=0.
\end{equation*}
Obviously, $\mathcal{E}_n$ and $E_n$ do not change.

\vspace{0.5 cm}

Set
\begin{equation}\label{eq:16***}
P_{\infty}=\left\{q^b\in (H^1(\Theta))^2:\, q^b\cdot \nu^b =0\,\, \text{on}\,\, \partial\Theta,\,\, q^b(0')=0 \right\},
\end{equation}
and
\begin{align}
E_{\infty}:q^b\in (H^1(\Theta))^2\mapsto \displaystyle 
&\nonumber \int_{\Theta}\left(\beta\left|\text{rot}\, q^b\right|^2+\left|\text{div}\, q^b\right|^2+\alpha\left(|q^b|^2-1\right)^2+|D\psi^b_{q^b}|^2\right)dx',\\
&\label{eq:18***}\quad +\int_{\Theta}\left(\int_{-1}^{0}(f^b_1, f^b_2) dx_3\cdot q^b\right)dx',
\end{align}
where $f^b_1$ and $f^b_2$ are defined in (\ref{eq:14}), and $\psi^{b}_{q^b}$ is the unique solution of 
\begin{equation}\label{eq:19***}
\begin{cases} 
      \psi^b_{p^b}\in H^1(\Theta) , \quad \displaystyle\int_{\Theta}\psi^b_{q^b}dx'=0, &   \\
      \displaystyle\int_{\Theta}\left(\left(-D_{x'}\psi^b_{q^b}+q^b\right)\cdot D_{x'} \psi^b\right)dx'=0, & 
 \end{cases}
\end{equation}
for $\psi^b\in H^1(\Theta)$.

\vspace{1.0 cm}

\begin{teo}\label{teo:3}
Assume (\ref{eq:1}) with $\ell=+\infty$ and $h^b_n\ll \sqrt{h^a_n}$, and (\ref{eq:14}). For every $n\in \NN$, let $(p^a_n, p^b_n)$ a solution of $(\ref{eq:13})$, and $\left(\phi^a_{(p^a_n, p^b_n)}, \phi^b_{(p^a_n, p^b_n)}\right)$ be the unique solution of $(\ref{eq:11})$ with $(p^a, p^b)=(p^a_n, p^b_n)$. Moreover, let $P_{\infty}$ and $E_{\infty}$ defined by (\ref{eq:16***}) and (\ref{eq:18***}), (\ref{eq:19***}) respectively. Then there exist an increasing sequence of positive integer numbers $\{n_i\}_{i\in\NN}$ and (in possible dependence on the subsequence) $(\hat{p}^b_1, \hat{p}^b_2)\in P_{\infty}$ such that
\begin{equation}\label{eq:20***}
\begin{array}{rcl}
          \displaystyle\frac{h^a_n}{\sqrt{h^b_n}}p^a_{n_i}\rightarrow 0 & \text{strongly in} & \left(H^1(\Omega^a)\right)^3\,  \text{ and strongly in}\,  \left(L^4(\Omega^a)\right)^3,\\
          p^b_{n_i}\rightarrow (\hat{p}^b_1,\hat{p}^b_2, 0) & \text{strongly in} & \left(H^1(\Omega^b)\right)^3\,  \text{ and strongly in}\,  \left(L^4(\Omega^b)\right)^3,
\end{array}
\end{equation}

\begin{equation}\label{eq:21***}
\begin{cases} 
      \displaystyle\left(\frac{1}{\sqrt{h^b_n}}\frac{\partial p^a_n}{\partial x_1}, \frac{1}{\sqrt{h^b_n}}\frac{\partial p^a_n}{\partial x_2}\right)\rightarrow(0,0) & \text{strongly in}\, \left(L^2(\Omega^a)\right)^{3}\times \left(L^2(\Omega^a)\right)^{3}, \\
      \displaystyle\frac{1}{h^b_n}\frac{\partial p^b_n}{\partial x_3}\rightarrow 0 & \text{strongly in}\, \left(L^2(\Omega^b)\right)^{3}, 
 \end{cases}
\end{equation}

\begin{equation}\label{eq:22***}
\begin{cases} 
      \displaystyle\left(\frac{h^a_n}{\sqrt{h^b_n}}\phi^a_{(p^a_{n_i},p^b_{n_i})}, \phi^b_{(p^a_{n_i},p^b_{n_i})} \right)\rightarrow(0, \psi^b_{\hat{p}^b}) & \text{strongly in}\, \left(H^1(\Omega^a)\right)\times \left(H^1(\Omega^b)\right), \\
     \displaystyle\left(\frac{1}{\sqrt{h^b_n}}\frac{\partial \phi^a_{(p^a_{n_i},p^b_{n_i})}}{\partial x_1} , \frac{1}{\sqrt{h^b_n}}\frac{\partial \phi^a_{(p^a_{n_i},p^b_{n_i})}}{\partial x_2} \right)\rightarrow(0,0) & \text{strongly in}\, \left(L^2(\Omega^a)\right)^2,\\
       \displaystyle\frac{1}{h^b_n}\frac{\partial  \phi^b_{(p^a_{n_i},p^b_{n_i})} }{\partial x_3}\rightarrow 0 & \text{strongly in}\, \left(L^2(\Omega^b)\right),  
 \end{cases}
\end{equation}
where $\hat{p}^b$ solves
\begin{equation}\label{eq:23***}
E_{\infty}(\hat{p}^b)=\min \{E_{\infty}(q^b):  q^b\in P_{\infty}\},
\end{equation}
 and $\psi^b_{\hat{p}^b}$ is the unique solution of (\ref{eq:19***}) with $q^b=\hat{p}^b$. Moreover
 \begin{equation}\label{eq:24***}
 \displaystyle\lim_n \frac{E_n((p^a_n, p^b_n))}{h^b_n}=E_{\infty}(\hat{p}^b).
\end{equation}
\end{teo}

\vspace{0.5 cm}

For brevity, we omit the proof of Theorem \ref{teo:3} which depevelops like that of Therem \ref{teo:2}. 

\section{The proof in the case $\ell\in]0, +\infty[$}\label{cha:4}
\vspace{1.0 cm}

\subsection{A priori estimates on polarization}
\vspace{0.5 cm}

\begin{prop}\label{prop:1}
Assume (\ref{eq:7}). For every $n\in\NN$, let $(p^a_n, p^b_n)$ be a solution of (\ref{eq:13}). Then, there exist a constant $c$ such that
\begin{align}
&\label{eq:25}||p^a_{n}||_{\left(L^4(\Omega^a)\right)^3}\leq c,\quad ||p^b_n||_{\left(L^4(\Omega^b)\right)^3}\leq c\quad \text{for all}\,\, n\in\NN,\\ 
 &\label{eq:26}||D^a_np^a_{n}||_{\left(L^2(\Omega^a)\right)^9}\leq c,\quad ||D^b_np^b_n||_{\left(L^2(\Omega^b)\right)^9}\leq c\quad \text{for all}\,\, n\in\NN.               
 \end{align}
 
 \vspace{0.5 cm}
 
 \dem
 Function $0$ belonging to $P_n$ gives
\begin{align}
&\nonumber\displaystyle\int_{\Omega^a}\left(\beta|\text{rot}_n^a\, p^a_n|^2+|\text{div}_n^a\, p^a_n|^2+\alpha\left(|p^a_n|^4-2|p^a_n|^2\right)+|D_n^a\phi^a_{(p^a_n, p^b_n)}|^2\right)dx\\
&\nonumber\quad +\frac{h^b_n}{(h^a_n)^2}\int_{\Omega^b}\left(\beta|\text{rot}_n^b\, p^b_2|^2+|\text{div}_n^b\, p^b_n|^2+\alpha\left(|p^b_n|^4-2|p^b_n|^2\right)^2+|D_n^b\phi^b_{(p^a, p^b)}|^2\right)dx,\\
&\label{eq:27}\quad\leq\displaystyle\frac{1}{2}\int_{\Omega^a}\left(|f^a_n|^2+|p^a_n|^2\right)dx+\frac{h^b_n}{(h^a_n)^2}\displaystyle\frac{1}{2}\int_{\Omega^b}\left(|f^b_n|^2+|p^b_n|^2\right)dx\quad \text{for all $n\in\NN$}.
\end{align}
Estimates (\ref{eq:27}) implies
\begin{align*}
&\nonumber\displaystyle\int_{\Omega^a}\alpha\left(|p^a_n|^4-\left(2+\frac{1}{2\alpha}\right)|p^a_n|^2\right)dx+\frac{h^b_n}{(h^a_n)^2}\int_{\Omega^b}\alpha\left(|p^b_n|^4-\left(2+\frac{1}{2\alpha}\right)|p^b_n|^2\right)dx\\
&\quad\leq\displaystyle\frac{1}{2}\int_{\Omega^a}\left(|f^a_n|^2+|p^a_n|^2\right)dx+\frac{h^b_n}{(h^a_n)^2}\displaystyle\frac{1}{2}\int_{\Omega^b}\left(|f^b_n|^2+|p^b_n|^2\right)dx\quad \text{for all $n\in\NN$},
\end{align*}
which gives
\begin{align}
&\nonumber\displaystyle\int_{\Omega^a}\alpha\left(|p^a_n|^2-\left(1+\frac{1}{4\alpha}\right)\right)^2dx+\frac{h^b_n}{(h^a_n)^2}\int_{\Omega^b}\alpha\left(|p^b_n|^2-\left(1+\frac{1}{4\alpha}\right)\right)^2dx\\
&\label{eq:28}\quad\leq\displaystyle\alpha\left(1+\frac{1}{4\alpha}\right)\left(|\Omega^a|+\frac{h^b_n}{(h^a_n)^2}|\Omega^b|\right)+\frac{1}{2}\int_{\Omega^a}|f^a_n|^2dx+\frac{h^b_n}{(h^a_n)^2}\displaystyle\frac{1}{2}\int_{\Omega^b}|f^b_n|^2dx\quad \text{for all $n\in\NN$}.
\end{align}
Then the estimate $(\ref{eq:25})$ follow from (\ref{eq:28}) , (\ref{eq:1}) with $\ell\in]0, +\infty[$ and (\ref{eq:14}). The estimate (\ref{eq:26}) follow from (\ref{eq:27}) with $\ell\in]0, +\infty[$, (\ref{eq:14}), (\ref{eq:25}), the continuous embedding of $L^4$ into $L^2$, and (\ref{eq:15}).
\lqqd
\end{prop}

\vspace{1.0 cm}

\begin{coro}\label{coro:1} 
Assume (\ref{eq:1}) with $\ell\in]0, +\infty[$ and (\ref{eq:7}). For every $n\in\NN$, let $(p^a_n, p^b_n)$ be a solution of (\ref{eq:13}). Let $P$ be defined in (\ref{eq:16}). Then there exist a subsequence of $\NN$, still denoted by $\{n\}$, and (in possible dependence on the subset) $(\hat{p}^a,\hat{p}^b)=(\hat{p}^a, (\hat{p}^b_1, \hat{p}^b_2))\in P$ such that
\begin{equation}\label{eq:29}
\begin{array}{rcl}
          p^a_{n}\rightharpoonup (0, 0, \hat{p}^a) & \text{weakly in} & \left(H^1(\Omega^a)\right)^3\,  \text{ and strongly in}\,  \left(L^4(\Omega^a)\right)^3,\\
          p^b_{n}\rightharpoonup (\hat{p}^b_1, \hat{p}^b_2, 0) & \text{weakly in} & \left(H^1(\Omega^b)\right)^3\,  \text{ and strongly in}\,  \left(L^4(\Omega^b)\right)^3.
\end{array}
\end{equation}
\end{coro}

\vspace{0.5 cm}

\dem Proposition \ref{prop:1} ensures that there exist a subsequence of $\NN$, still denoted by $\{n\}$, and (in possible dependence on the subsequence) $(\hat{p}^a_1, \hat{p}^a_2, \hat{p}^a_3)\in \left(H^1(\Omega^a)\right)^3$ independent of $x_1$ and $x_2$, and $(\hat{p}^b_1, \hat{p}^b_2, \hat{p}^b_3)\in \left(H^1(\Omega^b)\right)^3$ independent of  $x_3$ such that
\begin{equation}\label{eq:30}
\begin{array}{rcl}
          p^a_{n}\rightharpoonup (\hat{p}^a_1, \hat{p}^a_2, \hat{p}^a) & \text{weakly in} & \left(H^1(\Omega^a)\right)^3\,  \text{ and strongly in}\,  \left(L^4(\Omega^a)\right)^3,\\
          p^b_{n}\rightharpoonup (\hat{p}^b_1, \hat{p}^b_2, \hat{p}^b_3) & \text{weakly in} & \left(H^1(\Omega^b)\right)^3\,  \text{ and strongly in}\,  \left(L^4(\Omega^b)\right)^3,
\end{array}
\end{equation}
and $(\hat{p}^a_1, \hat{p}^a_2, \hat{p}^a_3)\cdot \nu^a=0$ on $\partial\Omega^a\setminus\left(\Theta\times\{0\}\right)$, $(\hat{p}^b_1, \hat{p}^b_2, \hat{p}^b_3)\cdot \nu^b=0$ on $\partial\Omega^b\setminus\left(\Theta\times\{0\}\right)$.

\vspace{0.5 cm}

In particular, this implies
\begin{equation*}
\begin{array}{rcl}
          \hat{p}^a_1=\hat{p}^a_2=0 & \text{in} & \Omega^a\,  \text{and}\,\,  \hat{p}^a(1)=0,\\
          \hat{p}^b_3=0 & \text{in} & \Omega^b\,  \text{and}\,\, (\hat{p}^b_1,\hat{p}^b_2)\cdot \nu^b =0\quad\text{on}\,\partial\Theta.
\end{array}
\end{equation*}
It remains to prove that
\begin{equation}\label{eq:31}
 \hat{p}^a(0)=(\hat{p}^b_1, \hat{p}^b_2)(0').
 \end{equation}
 
\vspace{0.5 cm}

For this end, we can argue as in \cite{REF18} as follows: Slim the proof of (\ref{eq:31}) in two claim. 

\vspace{0.5 cm}

\textit{Claim 1}: The existence of  three constants $c\in]0, +\infty[$, $\overline{x_3}\in]-1, 0[$ and of an increasing sequence of positive integers $\{i_k\}_{k\in\NN}$ such that
\begin{equation}\label{eq:32}
\displaystyle\int_{\Theta}\left|D_{x'}p^b_{n_{i_k}}(x', \overline{x_3})\right|^2 dx'\leq c,\quad \forall k\in\NN
\end{equation}
and
\begin{equation}\label{eq:33}
p^b_{n_{i_k}}(\cdot, \overline{x_3})\rightarrow \hat{p}^b\quad\text{strongly in}\, C^{0}\left(\Theta\right),
\end{equation}
as $k\to +\infty$.

\vspace{0.5 cm}
\dem
For every $i\in\NN$, set
\begin{align*}
&\rho_{i}: x_3\in]-1, 0[\mapsto \displaystyle\int_{\Theta}\left|p^b_{n_i}(x', x_3)\right|^2dx'\\
&\quad\displaystyle+\int_{\Theta}\left(\left|D_{x_1}p^b_{n_i}(x', x_3)\right|^2+\left|D_{x_2}p^b_{n_i}(x', x_3)\right|^2\right)dx'.
\end{align*}
Proposition \ref{prop:1}, (\ref{eq:26}) ensure there exist a subsequence of $\NN$, still denoted by $\{n_i\}$ and (in possible dependence on the subsequence). $(\xi^b_1, \xi^b_2, \xi^b_3)\in \left(L^2(\Omega^b)\right)^3$ such that 
\begin{equation}\label{eq:34}
\begin{cases} 
           \displaystyle D_{x_1}p^b_{n_ i}\rightharpoonup\xi^b_1,\quad  D_{x_2}p^b_{n_i}\rightharpoonup \xi^b_2 & \text{weakly in}\, L^2(\Omega^b),\\
       \displaystyle  \frac{1}{h^b_{n_i}}D_{x_3}p^b_{n_i}\rightharpoonup \xi^b_3 & \text{weakly in}\, L^2(\Omega^b).
 \end{cases}
\end{equation}
From Fatou's lemma and (\ref{eq:34}) and (\ref{eq:30}), it follows that
\begin{equation*}
\displaystyle\int_{-1}^{0}\lim_{i}\inf \rho_i(x_3)dx_3\leq \lim_{i}\inf \int_{-1}^{0}\rho_i (x_3)dx_3<+\infty.
\end{equation*} 
Consequently, there exist two constants $c\in]0,+\infty[$ and $\overline{x_3}\in]-1,0[$ and an insteresting sequence of positive integers $\{i_k\}_{k\in\NN}$ such that
\begin{equation*}
\rho_{i_k}(\overline{x_3})<c,\quad \forall k\in\NN
\end{equation*}
\textit{i.e.} estimate (\ref{eq:32}) holds true and, by virtue of the second convergence in (\ref{eq:30}), it follows that
\begin{equation}\label{eq:35}
p^b_{n_{i_k}}(\cdot, \overline{x_3})\rightharpoonup \hat{p}^b=(\hat{p}_1^b, \hat{p}_2^b) \quad \text{weakly in}\, (H^1\left(\Theta\right))^2
\end{equation}
as $k\to+\infty$, which provides (\ref{eq:33}).

\vspace{0.5 cm}

\textit{Claim 2}: 
\begin{equation}\label{eq:36}
\displaystyle\lim_{k}\int_{\Theta}p^b_{n_{i_k}}(h^a_{n_{i_k}}x', 0)dx'=|\Theta|p^b(0')=|\Theta|(p^b_1, p^b_2)(0').
\end{equation} 

\dem
We can split the integral in (\ref{eq:36}) in the following way:
\begin{align}
\displaystyle\int_{\Theta}p^b_{n_{i_k}}(h^a_{n_{i_k}}x', 0)dx'&\nonumber=\displaystyle\int_{\Theta}\left(p^b_{n_{i_k}}(h^a_{n_{i_k}}x', 0)-p^b_{n_{i_k}}(h^a_{n_{i_k}}x', \overline{x_3})\right)dx'\\
&\nonumber\quad+\displaystyle\int_{\Theta}\left(p^b_{n_{i_k}}(h^a_{n_{i_k}}x', \overline{x_3})-\hat{p}^b(h^a_{n_{i_k}}x'))\right)dx'\\
&\label{eq:37}\quad+\displaystyle\int_{\Theta}\hat{p}^b(h^a_{n_{i_k}}x')dx',\quad \forall k\in\NN,
\end{align}
and one will pass to the limit, as $k$ diverges, in each term of this descomposition. By using the last convergence in (\ref{eq:34}), there exist a constant $c\in]0, +\infty[$ such that
\begin{align}
&\nonumber\lim_{k}\sup \left|\displaystyle\int_{\Theta}\left(p^b_{n_{i_k}}(h^a_{n_{i_k}}x', 0)-p^b_{n_{i_k}}(h^a_{n_{i_k}}x', \overline{x_3})\right)dx'\right|\\
&\nonumber=\lim_{k} \sup\left|\displaystyle\int_{\Theta}\left(\int_{\overline{x_3}}^{0} D_{x_3}p^b_{n_{i_k}}(h^a_{n_{i_k}}x', x_3)dx_3\right)dx'\right|\\
&\nonumber\leq |\Omega^b|^{1/2}\lim_{k} \sup\left(\displaystyle\int_{\Omega^b}\left|D_{x_3}p^b_{n_{i_k}}(h^a_{n_{i_k}}x', x_3)\right|^2dx\right)^{1/2}\\
&\nonumber\leq |\Omega^b|^{1/2}\lim_{k} \sup\left(\displaystyle\frac{1}{(h^a_{n_{i_k}})^2}\int_{\Omega^b}\left|D_{x_3}p^b_{n_{i_k}}(x', x_3)\right|^2dx\right)^{1/2}\\
&\nonumber= |\Omega^b|^{1/2}\lim_{k} \sup \frac{h^b_{n_{i_k}}}{h^a_{n_{i_k}}} \lim_{k}\sup \left(\displaystyle\int_{\Omega^b}\left|\frac{1}{h^b_{n_{i_k}}}D_{x_3}p^b_{n_{i_k}}(x_1, x_2, x_3)\right|^2dx\right)^{1/2}\\
&\label{eq:38}\leq |\Omega^b|^{1/2}c\lim_{k}(h^b_{n_{i_k}})^{1/2} \lim_{k}\left(\frac{h^b_{n_{i_k}}}{(h^a_{n_{i_k}})^2}\right)^{1/2}=0.
\end{align}
The convergence in (\ref{eq:33}) implies that
\begin{align}
&\nonumber\displaystyle\lim_{k}\sup\left|\displaystyle\int_{\Theta}\left(p^b_{n_{i_k}}(h^a_{n_{i_k}}x', \overline{x_3})-\hat{p}^b_{1}(h^a_{n_{i_k}}x')\right)dx'\right|\\
&\nonumber \leq \displaystyle\lim_{k}\sup\displaystyle\int_{\Theta}\left|\left(p^b_{n_{i_k}}(h^a_{n_{i_k}}x', \overline{x_3})-\hat{p}^b_{1}(h^a_{n_{i_k}}x')\right)\right|dx'\\
&\nonumber = \displaystyle\lim_{k}\sup\left(\displaystyle\frac{1}{(h^a_{n_{i_k}})^2}\int_{h^a_{n_{i_k}}\Theta}\left|\left(p^b_{n_{i_k}}(x', \overline{x_3})-\hat{p}^b(x')\right)\right|dx'\right)\\
&\label{eq:40}\leq\lim_{k}||p^b_{n_{i_k}}(\cdot, \overline{x_3})-\hat{p}^b(\cdot)||_{L^{\infty}(\Theta)}=0.
\end{align}
Since $\hat{p}^b\in C^0(\Theta)$, it follows that
\begin{equation}\label{eq:41}
\displaystyle\lim_{k}\int_{\Theta}\hat{p}^b(h^a_{n_{i_k}}x')dx'=|\Theta|\hat{p}^{b}(0').
\end{equation}
By passing to the limit in (\ref{eq:37}), as $k$ diverges, and taking into account $(\ref{eq:38})-(\ref{eq:41})$, one obtains (\ref{eq:36}).

\vspace{0.5 cm}

Finally, the junction condition (\ref{eq:31}) is obtained passing to the limit, as $k$ diverges, in
\begin{equation*}
 \displaystyle\int_{\Theta}p^a_{n}(x', 0)dx'=\int_{\Theta}p^b_{n}(h^a_nx', 0)dx',
\end{equation*}
and using the first convergence en (\ref{eq:30}) and $(\ref{eq:36})$.

\lqqd

\vspace{0.5 cm}
\subsection{A convergence result for problem (\ref{eq:11})}
\vspace{0.5 cm}

In this subsection we study the asymptotic behavior of the non local term generated by the potential solution of problem (\ref{eq:11}). Let 
\begin{equation}
U_{\text{reg}}\label{eq:42}=\left\{(\psi^a, \psi^b)\in C^1([0,1])\times C(\overline{\Theta}): \psi^b|_{\overline{\Theta}}\in C^1(\overline{\Theta}),\quad \psi^a(0)=\psi^b(0') \right\}.
\end{equation}
\vspace{0.5 cm}

\begin{prop}\label{prop:2} Let $U$ and $U_{\text{reg}}$ be defined in (\ref{eq:17}) and (\ref{eq:42}) respectively. Then $U_{\text{reg}}$ is dense in $U$.
\end{prop}
\vspace{0.5 cm}

\dem
Let $(\psi^a, \psi^b)\in U$. The goal is to find a sequence $\{(\psi^a_n, \psi^b_n)\}_{n\in\NN}\subset U_{\text{reg}}$ such that
\begin{equation}\label{eq:43}
(\psi^a_n, \psi^b_n)\rightarrow (\psi^a, \psi^b)\quad\text{strongly in}\, H^1(]0,1[)\times H^1(\Theta).
\end{equation} 
To this end, split $\psi^b=\psi^{e}+\psi^{o}$ in the even part and in the odd part with respect to $x_1$. Note that $\psi^{e}$ and $\psi^{o}$ belong to $H^1(\Theta)$, and
\begin{equation*}
\psi^{e}(0')=\psi^{b}(0')=\psi^{a}(0), \quad \psi^{o}(0')=0.
\end{equation*}
Consequently, a convolution argument allows us to build three sequences $\{\zeta^a_n\}_{n\in\NN}\subset C^{\infty}(]0,1[)$, $\{\zeta^e_n\}_{n\in\NN}\subset C(\Theta)$ and $\{\zeta^o_n\}_{n\in\NN}\subset C^{\infty}(\Theta)$ such that
\begin{equation}\label{eq:44}
\begin{cases} 
          \left\{\zeta^e_{n}|_{\overline{\Theta}}\right\}_{n\in\NN}\subset C^{\infty}(\overline{\Theta}), &   \\
\zeta^a_n\rightarrow \psi^a\quad\text{strongly in}\, H^1(]0, 1[),\\
\zeta^e_n\rightarrow \psi^e\quad\text{strongly in}\, H^1(\Theta),\quad \zeta^0_{n}\rightarrow \zeta^0\quad \text{strongly in}\, H^1(\Theta),\\
\zeta^a_n(0)=\zeta^o_n(0'),\quad \zeta^{o}_{n}(0')=0\quad\text{for all}\,n\in\NN.
   \end{cases}
\end{equation}
This implies (\ref{eq:43}), setting $\psi^a_n=\zeta^a_n$ and $\psi^b_n=\zeta^e_n+\zeta^o_n$.
\lqqd
\vspace{1.0 cm}

\begin{prop}\label{prop:3}
Assume (\ref{eq:1}) with $\ell\in]0, +\infty[$. Let $\{(q^a_n, q^b_n)\}_{n\in\NN}\subset \left(L^2(\Omega^a)\right)^3\times (L^2(\Omega^b))^3$, and let $(q^a,q^b)\in \left(L^2(\Omega^a)\right)\times (L^2(\Omega^b))^2$ such that $q^a$ is independent of $(x_1,x_2)$, $q^b=(q^b_1, q^b_2)$ is independent of $x_3$ and 
\begin{equation}\label{eq:45}
(q^a_n, q^b_n)\rightarrow ((0, 0, q^a),(q^b_1, q^b_2, 0))\quad\text{strongly in}\, \left(L^2(\Omega^a)\right)^3\times (L^2(\Omega^b))^3. 
\end{equation}
Moreover, for $n\in\NN$ let $(\phi^a_{(q^a_n, q^b_n)}, \phi^b_{(q^a_n, q^b_n)} )$ be the unique solution of 
\begin{equation}\label{eq:46}
\begin{cases} 
      \left(\phi^a_{(q^a_n, q^b_n)},\phi^b_{(q^a_n, q^b_n)} \right)\in U_n, \quad \displaystyle\int_{\Omega^a}\phi^a_{(q^a_n, q^b_n)}dx=0, &   \\
      \displaystyle\int_{\Omega^a}\left(\left(-D^a_n\phi^a_{(q^a_n, q^b_n)}+q^a_n\right)\cdot D^a_n \phi^a\right)dx+\frac{h^b_n}{(h^a_n)^2}\displaystyle\int_{\Omega^b}\left(\left(-D^b_n\phi^b_{(q^a_n, q^b_n)}+q^b_n\right)\cdot D^b_n \phi^b\right)dx=0, & 
         \end{cases}
\end{equation}
$(\phi^a, \phi^b)\in U_n,$ where $U_n$ is defined in (\ref{eq:9}). Furthermore, let $(\psi^a_{(q^a, q^b)}, \psi^b_{(q^a, q^b)})$ be the unique solution of (\ref{eq:19}). Then one has
\begin{equation}\label{eq:47}
\begin{cases} 
      \displaystyle\left(\phi^a_{(q^a_{n},q^b_{n})}, \phi^b_{(q^a_{n}, q^b_{n})} \right)\rightarrow(\psi^a_{(q^a, q^b)}, \psi^b_{(q^a, q^b)}) & \text{strongly in}\, \left(H^1(\Omega^a)\right)\times \left(H^1(\Omega^b)\right), \\
     \displaystyle\left(\frac{1}{h^a_n}\frac{\partial \phi^a_{(q^a_{n}, q^b_{n})}}{\partial x_1}, \frac{1}{h^a_n}\frac{\partial \phi^a_{(q^a_{n}, q^b_{n})}}{\partial x_2} \right)\rightarrow(0,0) & \text{strongly in}\, \left(L^2(\Omega^a)\right)^2,\\
       \displaystyle \frac{1}{h^b_n}\frac{\partial  \phi^b_{(q^a_{n}, q^b_{n})} }{\partial x_3}\rightarrow 0 & \text{strongly in}\,  L^2(\Omega^b).  
 \end{cases}
\end{equation}
\end{prop}

\vspace{0.5 cm}

\dem
In this proof, $c$ denotes any positive constant independent of $n\in\NN$. Choosing $(\phi^a, \phi^b)=\left(\phi^a_{(q^a_{n}, q^b_{n})}, \phi^b_{(q^a_{n}, q^b_{n})} \right)$ as test function in (\ref{eq:46}) applying Young inequality, and using (\ref{eq:45}) give
\begin{equation}\label{eq:48}
||D^a_n\phi^a_{(q^a_{n}, q^b_{n})}||_{(L^2(\Omega^a))^3}\leq c,\quad ||D^b_n\phi^b_{(q^a_{n}, q^b_{n})}||_{(L^2(\Omega^b))^3}\leq c\quad \text{for all}\, n\in\NN.
\end{equation}
the first estimate in (\ref{eq:48}) implies
\begin{equation}\label{eq:49}
 ||\phi^a_{(q^a_{n}, q^b_{n})}||_{H^1(\Omega^a)}\leq c\quad \text{for all}\, n\in\NN,
\end{equation}
since $\displaystyle\int_{\Omega^a}\phi^a_{(q^a_{n}, q^b_{n})}dx=0$. The next step is devoted to proving
\begin{equation}\label{eq:50}
 ||\phi^b_{(q^a_{n}, q^b_{n})}||_{H^1(\Omega^b)}\leq c\quad \text{for all}\, n\in\NN.
\end{equation}
The junction condition in (\ref{eq:9}) gives  
\begin{align}
\displaystyle\int_{h^a_n\Theta}\left|\phi^b_{(q^a_{n}, q^b_{n})}(x', 0)\right|^2 dx_1dx_2&\nonumber=(h^a_n)^2\int_{\Theta}\left|\phi^b_{(q^a_{n}, q^b_{n})}( h^a_n x', 0)\right|^2 dx'\\
&\label{eq:51}=(h^a_n)^2\int_{\Theta}\left|\phi^a_{(q^a_{n}, q^b_{n})}(x_1, x_2, 0)\right|^2 dx',\,\text{for all}\, n\in\NN.
\end{align}
Then, (\ref{eq:51}) and (\ref{eq:49}) provide
\begin{equation*}
||\phi^b_{(q^a_{n}, q^b_{n})}||_{L^2\left(h^a_n\Theta\times\{0\}\right)}\leq h^a_n c,\quad\text{for all}\,n\in\NN,
\end{equation*}
which implies
\begin{equation}\label{eq:52}
||\phi^b_{(q^a_{n}, q^b_{n})}||_{L^2\left(h^a_n\Theta\times]-1, 0[\right)}\leq h^a_nc,\quad\text{for all}\,n\in\NN,
\end{equation}
by virtue of the (\ref{eq:48}) ensures
\begin{equation*}
||\phi^b_{q_n}||_{L^2\left(\{0\}\times\{0\}\times]-1,0[\right)}\leq c,\quad\text{for all}\,n\in\NN,
\end{equation*}
which combined again with the second estimates in (\ref{eq:48}) proves (\ref{eq:50}). The estimates (\ref{eq:48}), (\ref{eq:49}), and (\ref{eq:50}) ensure the existence of a subsequence of $\NN$, still denotes by $\{n\}$ and (in possible dependence on the subsequence) $(\tau^a, \tau^b)\in U$, $(\xi^a, \zeta^a)\in (L^2(\Omega^a))^2$ and $\zeta^b\in L^2(\Omega^b)$ such that
\begin{align} 
      \displaystyle\left(\phi^a_{(q^a_{n},q^b_{n})}, \phi^b_{(q^a_{n}, q^b_{n})} \right)\rightharpoonup(\tau^a, \tau^b) &\label{eq:53}\, \text{weakly in}\, \left(H^1(\Omega^a)\right)\times \left(H^1(\Omega^b)\right), \\
     \displaystyle\left(\frac{1}{h^a_n}\frac{\partial \phi^a_{(q^a_{n}, q^b_{n})}}{\partial x_1}, \frac{1}{h^a_n}\frac{\partial \phi^a_{(q^a_{n}, q^b_{n})}}{\partial x_2} \right)\rightharpoonup (\xi^a,\zeta^a) &\label{eq:54}\, \text{weakly in}\, \left(L^2(\Omega^a)\right)^2,\\
       \displaystyle  \frac{1}{h^b_n}\frac{\partial  \phi^b_{(q^a_{n}, q^b_{n})} }{\partial x_3}\rightharpoonup\zeta^b &\label{eq:55}\, \text{weakly in}\,  L^2(\Omega^b),  
 \end{align}
and
\begin{equation}\label{eq:56}
\displaystyle\int_{0}^{1}\tau^a dx_3=0.
\end{equation}
Note that the junction condition $\tau^a(0)=\tau^b(0')$ can be obtained arguing as the proof of (\ref{eq:31}). 

\vspace{0.5 cm}

The next step is devoted to proving that $(\tau^a, \tau^b)$ solves (\ref{eq:19}). To this end, for every couple $(\psi^a, \psi^b)\in U_{\text{reg}}$ where $U_{\text{reg}}$ is defined in (\ref{eq:42}), consider a sequence $\{\mu_n\}_{n\in\NN}\subset H^1(\Omega^a)$ (depending on $(\psi^a, \psi^b)$) such that
\begin{equation}\label{eq:57}
\begin{cases} 
      \displaystyle\left(\mu_n, \psi^b \right)\in U_n & \text{for all}\, n\in\NN, \\
     \displaystyle \mu_n \rightarrow \psi^a & \text{strongly in}\, L^2(\Omega^a),\\
       \displaystyle\left(\frac{1}{h^a_n}D_{x_1}\mu_n, \frac{1}{h^a_n}D_{x_2}\mu_n, D_{x_3}\mu_n\right)\rightarrow(0,0, D_{x_3}\psi^a) & \text{strongly in}\, \left(L^2(\Omega^a)\right)^3.  
 \end{cases}
\end{equation}
For instance, setting
\begin{equation*}
\mu_n(x)=
\begin{cases} 
      \displaystyle\psi^a(x_3) & \text{if}\, x=(x', x_3)\in\Theta\times]h^a_n,1[, \\
       \displaystyle\psi^a(h^a_n)\frac{x_3}{h^a_n}+\psi^b(h^a_n x')\frac{h^a_n-x_3}{h^a_n} & \text{if}\, x=(x', x_3)\in\Theta\times[0, h^a_n],  
 \end{cases}
\end{equation*}
the first two properties in (\ref{eq:57}) can be immediately verified by the properties of $U_{\text{reg}}$, while the last ones follows from
\begin{align*}
\displaystyle\int_{\Theta\times]0, h^a_n[} \left|\frac{1}{h^a_n}D_{x_1}\mu_n\right|^2dx&=\int_{\Theta\times]0, h^a_n[}\left|D_{x_1}\psi^b(h^a_n x')\left(1-\frac{x_3}{h^a_n}\right)\right|^2dx\\
&\displaystyle=\int_{\Theta}\left|D_{x_1}\psi^b(h^a_n x')\right|^2dx'\int_{]0, h^a_n[}\left(1-\frac{x_3}{h^a_n}\right)^2dx_3\\
&\leq ||\psi^b||^2_{W^{1,\infty}(\Theta)}h^a_n,\quad\text{for all}\, n\in\NN.
\end{align*}

\begin{align*}
\displaystyle\int_{\Theta\times]0, h^a_n[} \left|\frac{1}{h^a_n}D_{x_2}\mu_n\right|^2dx&=\int_{\Theta\times]0, h^a_n[}\left|D_{x_2}\psi^b(h^a_n x')\left(1-\frac{x_3}{h^a_n}\right)\right|^2dx\\
&\displaystyle=\int_{\Theta}\left|D_{x_2}\psi^b(h^a_n x')\right|^2dx'\int_{]0, h^a_n[}\left(1-\frac{x_3}{h^a_n}\right)^2dx_3\\
&\leq ||\psi^b||^2_{W^{1,\infty}(\Theta)}h^a_n,\quad\text{for all}\, n\in\NN.
\end{align*}

\begin{align*}
\displaystyle\int_{\Theta\times]0, h^a_n[} \left| D_{x_3}\mu_n\right|^2dx&=\int_{\Theta\times]0, h^a_n[}\left|\psi^a(h^a_n)\frac{1}{h^a_n}-\psi^{b}(h^a_nx')\frac{1}{h^a_n}\right|^2dx\\
&\displaystyle=\int_{\Theta}\frac{1}{h^a_n}\left|\psi^a(h^a_n)-\psi^{b}(h^a_nx')\right|^2dx'\\
&\displaystyle=\int_{\Theta}\frac{1}{h^a_n}\left|\psi^a(h^a_n)-\psi^a(0)+\psi^b(0')-\psi^{b}(h^a_nx')\right|^2dx'\\
&\leq2\left(|\Theta|\left|\frac{d\psi^a}{dx_3}|_{x_3=0}\right|^2+||\psi^b||^2_{W^{1,\infty}(\Theta)} \right)h^a_n,\quad\text{for all}\,n\in\NN,
\end{align*}
where again the properties of $U_{\text{reg}}$ played a crucial role.

\vspace{0.5 cm}

Now, fixing $(\psi^a, \psi^b)\in U_{\text{reg}}$, choosing $(\mu_n, \psi^b)$ as test function in $(\ref{eq:46})$ with $\mu_n$ satisfying (\ref{eq:57}), pasing to the limit as $n$ diverges, and using (\ref{eq:45}), (\ref{eq:53})-(\ref{eq:55}), one obtains
\begin{equation*}
 \displaystyle\int_{\Omega^a}\left(\left(-D^a_n\phi^a_{(q^a_n, q^b_n)}+q^a_n\right)\cdot D^a_n \mu_n\right)dx+\displaystyle\frac{h^b_n}{(h^a_n)^2}\int_{\Omega^b}\left(\left(-D^b_n\phi^b_{(q^a_n, q^b_n)}+q^b_n\right)\cdot D^b_n \psi^b\right)dx=0,
\end{equation*}
taking the limit as $n$ diverges
\begin{equation*}
 \displaystyle\int_{\Theta\times]0, 1[}\left(\left(-D_{x_3}\tau^a+q^a\right)\cdot D_{x_3}\psi^a\right)dx_1 dx_2 dx_3+\ell\displaystyle\int_{\Theta\times]-1, 0[}\left(\left(-D_{x'}\tau^b+q^b\right)\cdot D_{x'} \psi^b\right)dx_1dx_2dx_3=0,
\end{equation*}
\begin{equation}\label{eq:58}
 \displaystyle|\Theta|\int_{0}^{1}\left(\left(-D_{x_3}\tau^a+q^a\right)\cdot D_{x_3}\psi^a\right)dx_3+\ell\displaystyle\int_{\Theta}\left(\left(-D_{x'}\tau^b+q^b\right)\cdot D_{x'} \psi^b\right)dx'=0.
\end{equation}
By virtue of Proposition \ref{prop:2}, equation (\ref{eq:58}) holds true also with any test function in $U$. Consequently, thanks to (\ref{eq:56}), $(\tau^a, \tau^b)$ is the unique solution of (\ref{eq:19}), $\textit{i.e.}$
\begin{equation}\label{eq:59}
(\tau^a, \tau^b)=\left(\psi^a_{(q^a, q^b)}, \psi^b_{(q^a, q^b)}\right)\quad\text{a.e. in}\,\Omega^a\times\Omega^b.
\end{equation}
Finally, using (\ref{eq:46}), (\ref{eq:45}), (\ref{eq:59}), (\ref{eq:54})-(\ref{eq:55}), (\ref{eq:52}) and an l.s.c. argument, one has
\begin{align*}
&\displaystyle\lim_{n}\left(\int_{\Omega^a}\left|D^a_n\phi^a_{(q^a_n, q^b_n)}\right|^2dx+\frac{h^b_n}{(h^a_n)^2}\int_{\Omega^b}\left|D^a_n\phi^b_{(q^a_n, q^b_n)}\right|^2dx\right)\\
&\quad\geq \int_{\Omega^a}\left|\xi^a\right|^2dx+\int_{\Omega^a}\left|D_{x_3}\tau^a\right|^2dx+\int_{\Omega^a}\left|\zeta^a\right|^2dx+\ell\int_{\Omega^b}\left|D_{x'}\tau^b\right|^2dx+\ell\int_{\Omega^b}\left|\zeta^b\right|^2dx\\
&\quad= \int_{\Omega^a}\left|\xi^a\right|^2dx+|\Theta|\int_{0}^{1}\left|D_{x_3}\tau^a\right|^2dx_3+\int_{\Omega^a}\left|\zeta^a\right|^2dx+\ell\int_{\Theta}\left|D_{x'}\tau^b\right|^2dx'+\ell\int_{\Omega^b}\left|\zeta^b\right|^2dx.
\end{align*}
In the other hand,
\begin{align*}
&\displaystyle\lim_{n}\left(\int_{\Omega^a}\left|D^a_n\phi^a_{(q^a_n, q^b_n)}\right|^2dx+\frac{h^b_n}{(h^a_n)^2}\int_{\Omega^b}\left|D^a_n\phi^b_{(q^a_n, q^b_n)}\right|^2dx\right)\\
&\underset{(\ref{eq:46})}{=}\displaystyle\lim_{n}\left(\int_{\Omega^a}D^a_n\phi^a_{(q^a_n, q^b_n)}\cdot q^a_ndx+\frac{h^b_n}{(h^a_n)^2}\int_{\Omega^a}D^b_n\phi^b_{(q^a_n, q^b_n)}\cdot q^b_ndx\right)\\
&=|\Theta|\displaystyle\int_{0}^{1}D_{x_3}\psi^a_{(q^a, q^b)}\cdot q^adx_3+\ell\int_{\Theta}D_{x'}\psi^b_{(q^a, q^b)}\cdot q^bdx'\\
&\underset{(\ref{eq:19})}{=}|\Theta|\displaystyle\int_{0}^{1}\left|D_{x_3}\psi^a_{(q^a, q^b)}\right|^2dx_3+\ell\int_{\Theta}\left|D_{x'}\psi^b_{(q^a, q^b)}\right|^2dx',\\
&=\displaystyle|\Theta|\int_{0}^{1}\left|D_{x_3}\tau^a\right|^2dx_3+\ell\int_{\Theta}\left|D_{x'}\tau^b\right|^2dx',
\end{align*}
which implies that
\begin{equation*}
\xi^a=\zeta^a=\zeta^b=0,
\end{equation*}
and the convergence (\ref{eq:53})-(\ref{eq:55}) are strong. Note that previous convergences hold true for the whole sequence, since the limits are uniquely identified.
\lqqd

\vspace{1.0 cm}
\subsection{The proof of Theorem \ref{teo:1}}
\vspace{0.5 cm}
Let 
\begin{equation}
P_{\text{reg}}\label{eq:60}=C^1_0\left(]0,1[\right)\times \left(C^1_0(\Theta\setminus\{0'\})\right)^2.
\end{equation}

Then the following result is obvious.

\begin{prop}\label{prop:5}
Let $P$ and $P_{\text{reg}}$ defined in (\ref{eq:16}) and (\ref{eq:60}) respectively. Then $P_{\text{reg}}$ is dense in $P$.
\end{prop}

\subsubsection*{Proof of Theorem \ref{teo:1}} In what follows, $p^a_{n,i}$ (resp. $p^b_{n,i}$) denotes the $i-$th component, $i=1, 2, 3$, of $p^a_n$ (resp. $p^b_n$). Corollary \ref{coro:1} asserts that there exist a subsequence of $\NN$, still denoted by $\{n\}$, and (in possible dependence on the subsequence) $(\hat{p}^a, \hat{p}^b)\in P$ satisfying (\ref{eq:29}). The next step is devoted to proving the existence of a subsequence of $\NN$, still denoted by $\{n\}$, and (in possible dependence on the subsequence) of
\begin{align*}
&z^a\in L^2(H^1(\Theta), ]0,1[)\times (L^2(H^1_m(\Theta), ]0,1[))^2\\
&z^b\in(L^2(\Theta, H^1_m(]-1, 0[)))^2\times L^2(\Theta, H^1(]-1, 0[))
\end{align*}
\begin{equation}\label{eq:61}
\begin{cases}
z^{a}_1(\partial\Theta, \cdot)=0,\quad \text{on}\,]0,1[ &\\
z^{a}_2(\partial\Theta, \cdot)=0,\quad \text{on}\,]0,1[ &\\
z^{b}_2(\partial\Theta, \cdot)=0,\quad \text{on}\,]-1,0[ &\\
z^{b}_3(\cdot, 0)=0=z^b_{3}(\cdot, -1),\quad \text{on}\,\Theta&
\end{cases}
\end{equation}
 such that
\begin{equation}\label{eq:62}
\begin{cases} 
      &  \displaystyle\frac{1}{h^a_n}D_{(x_1,x_2)} p^a_n \rightharpoonup D_{(x_1, x_2)} z^a\,\, \text{weakly in}\, (L^2(\Omega^a))^6,\\
       &  \displaystyle\frac{1}{h^b_n}D_{x_3} p^b_n \rightharpoonup D_{x_3} z^b\,\, \text{weakly in}\, (L^2(\Omega^b))^3.
        \end{cases}
\end{equation}
Indeed, the boundary condition on $p^a_n$ and the Poincar\'e inequality give, for $i=1,2$,
\begin{equation}\label{eq:63}
\displaystyle\left|\left|\frac{1}{h^a_n}p^a_{n,i}(\cdot, \cdot, x_3)\right|\right|_{H^1(\Theta)}\leq \frac{c}{h^a_n}\left|\left|\left( \frac{\partial p^a_{n,i}(\cdot, \cdot, x_3) }{\partial x_1}, \frac{\partial p^a_{n,i}(\cdot, \cdot, x_3)}{\partial x_2}\right)\right|\right|_{L^2(\Theta)}
\end{equation}
for a.e. $x_3\in]0,1[$ and all $n\in\NN$, where $c$ is the Poincar\'e constant in $H^1_0(\Theta)$, while, for $i=3$, the Poincar\'e- Wirtinger inequality gives
\begin{align}
&\nonumber\displaystyle\left|\left|\frac{1}{h^a_n}\left(p^a_{n,3}(\cdot, \cdot, x_3)-\int_{\Theta}p^a_{n, 3}(x_1, x_2, x_3)dx_1dx_2\right)\right|\right|_{H^1_m(\Theta)}\\
&\label{eq:64}\displaystyle\quad\leq \frac{c'}{h^a_n}\left|\left|\left( \frac{\partial p^a_{n,3}(\cdot, \cdot, x_3) }{\partial x_1}, \frac{\partial p^a_{n,3}(\cdot, \cdot, x_3)}{\partial x_2}\right)\right|\right|_{L^2(\Theta)}
\end{align}
for a.e. $x_3\in]0,1[$, and all $n\in\NN$, where the subscript ``$m$'' means zero average, and $c'$ is the Poincar\'e-Wirtinger constant in $H^1_{m}(\Theta)$. Integrating (\ref{eq:63}) and (\ref{eq:64}) over $x_3\in]0, 1[$ and using the first estimate in (\ref{eq:26}) imply the first convergence in (\ref{eq:62}).Similarly, one proves the remaining convergences in $(\ref{eq:62})$. 

The next step is devoted to identifying $\hat{p}^a, \hat{p}^b, z^a$, and $z^b$. To this end, for every
\begin{equation*}
(q^a, q^b)\in P_{\text{reg}},
\end{equation*}
where $P_{\text{reg}}$ is defined in (\ref{eq:60}), consider a sequence $\{v_n\}_{n\in\NN}\subset (H^1(\Omega^a))^3$ (depending on $(q^a, q^b)$) such that
\begin{equation}\label{eq:65}
\begin{cases} 
      \displaystyle\left( v_n, (q^b_1, q^b_2, 0) \right)\in P_n & \text{for all}\quad n\in\NN, \\
     \displaystyle v_n \rightarrow (0, 0, q^a) & \text{strongly in}\, (L^4(\Omega^a))^3,\\
       \displaystyle\left(\frac{1}{h^a_n}D_{x_1} v_n, \frac{1}{h^a_n}D_{x_2} v_n, D_{x_3} v_n\right)\rightarrow(0,0, D_{x_3} q^a) & \text{strongly in}\, \left(L^2(\Omega^a)\right)^9.  
 \end{cases}
\end{equation}
For instance, setting
\begin{align*}
v_n(x)&=(0, 0, q^a)(x_3)\quad\text{if}\, x=(x', x_3)\in \Theta\times]h^a_n, 1[\,\text{and}\\
v_n(x)&=\left(q^b(h^a_nx')\frac{h^a_n-x_3}{h^a_n}, 0, q^a(x_3)\right)\quad\text{if}\, x=(x', x_3)\in \Theta\times[0, h^a_n],
\end{align*}
the first two properties in (\ref{eq:65}) can be immediately verified by virtue of the properties of $P_{\text{reg}}$, while the last ones follows from
\begin{align*}
\displaystyle\int_{\Theta\times]0, h^a_n[} \left|\frac{1}{h^a_n}D_{x_1}v_n\right|^2dx&=\int_{\Theta\times]0, h^a_n[}\left|D_{x_1} q^b(h^a_n x')\left(\frac{h^a_n-x_3}{h^a_n}\right)\right|^2dx\\
&\displaystyle=\int_{\Theta}\left|D_{x_1} q^b(h^a_n x')\right|^2dx'\int_{]0, h_n[}\left(1-\frac{x_3}{h^a_n}\right)^2dx_3\\
&\leq ||q^b||^2_{W^{1,\infty}(\Theta)}h^a_n,\quad\text{for all}\, n\in\NN.
\end{align*}
\begin{align*}
\displaystyle\int_{\Theta\times]0, h^a_n[} \left|\frac{1}{h^a_n}D_{x_2}v_n\right|^2dx&=\int_{\Theta\times]0, h^a_n[}\left|D_{x_2} q^b(h^a_n x')\left(\frac{h^a_n-x_3}{h^a_n}\right)\right|^2dx\\
&\displaystyle=\int_{\Theta}\left|D_{x_2} q^b(h^a_n x')\right|^2dx'\int_{]0, h_n[}\left(1-\frac{x_3}{h^a_n}\right)^2dx_3\\
&\leq ||q^b||^2_{W^{1,\infty}(\Theta)}h^a_n,\quad\text{for all}\, n\in\NN.
\end{align*}
\begin{align*}
\displaystyle\int_{\Theta\times]0, h^a_n[} \left| D_{x_3} v_n\right|^2dx&=\int_{\Theta\times]0, h^a_n[}\left|q^b(h^a_nx')\frac{1}{h^a_n}\right|^2dx+\int_{\Theta\times]0, h^a_n[} \left|D_{x_3}q^{a}(x_3)\right|^2dx\\
&\leq\int_{\Theta\times]0, h^a_n[}\left|\frac{1}{h^a_n}\left(q^b(h^a_nx')-q^b(0')\right)\right|^2dx\\
&\quad +\int_{\Theta\times]0, h^a_n[} \left|D_{x_3}q^{a}(x_3)\right|^2dx\\
&\leq\left(||q^b||^2_{W^{1,\infty}(]-\frac{1}{2},\frac{1}{2}[)}+||q^a||^2_{W^{1,\infty}(]0, 1[)}\right)h^a_n,\quad\text{for all}\,n\in\NN,
\end{align*}
where again the properties of $P_{\text{reg}}$ played a crucial role.

\vspace{0.5 cm}

Now, fixing $(q^a, q^b)=(q^a, (q^b_1, q^b_2))\in P_{\text{reg}}$, and choosing $(v_n, (q^b_1, q^b_2, 0))$ as test function in (\ref{eq:13}) with $v_n$ satisfying (\ref{eq:65}) give
\begin{equation}\label{eq:66}
\displaystyle\frac{1}{(h^a_n)^2}E_n\left((p^a_n, p^b_n)\right)\leq \frac{1}{(h^a_n)^2}E_n\left(v_n, (q^b_1, q^b_2, 0)\right)\quad\text{for all}\,n\in\NN.
\end{equation}
Then, passing to the limit in (\ref{eq:66}), as $n$ diverges, and using (\ref{eq:14}), (\ref{eq:29}), (\ref{eq:62}), (\ref{eq:65}), Proposition \ref{prop:3} and a l.s.c. argument imply

\begin{align}
&\nonumber\displaystyle\int_{\Omega^a}\left(\beta\left(\left|\frac{\partial z^a_3}{\partial x_2}\right|^2+\left|\frac{\partial z^a_3}{\partial x_1}\right|^2+\left|\frac{\partial z^a_2}{\partial x_1}-\frac{\partial z^a_1}{\partial x_2}\right|^2\right)+\left(\left|\frac{\partial z^a_1}{\partial x_1}+\frac{\partial z^a_2}{\partial x_2}+\text{div}\,\hat{p}^a\right|^2\right)\right)dx\\
&\nonumber\quad+|\Theta|\displaystyle\int_{0}^{1}\left(\alpha(|\hat{p}^a|^2-1)^2+|D_{x_3}\psi^a_{(\hat{p}^a, \hat{p}^b)}|^2\right)dx_3+\int_{0}^{1}\left(\int_{\Theta}f^a_3 dx_1dx_2\cdot \hat{p}^a \right)dx_3\\
&\nonumber\quad+\ell\displaystyle\int_{\Omega^b}\left(\beta\left(\left|\frac{\partial z^b_2}{\partial x_3}\right|^2+\left|\frac{\partial z^b_1}{\partial x_3}\right|^2\right)+ \left(\left|\frac{\partial z^b_3}{\partial x_3}+\text{div}\,\hat{p}^b\right|^2\right)\right)dx\\
&\nonumber\quad+\ell\displaystyle\int_{\Theta}\left(\alpha(|\hat{p}^b|^2-1)^2+|D_{x'}\psi^b_{(\hat{p}^a, \hat{p}^b)}|^2\right)dx'+\ell\int_{\Theta}\left(\int_{-1}^{0}(f^b_1, f^b_2) dx_3\cdot \hat{p}^b \right)dx'\\
&\nonumber\leq\lim_{n}\inf\frac{E_n\left((p^a_n, p^b_n)\right)}{(h^a_n)^2}
\label{eq:67}\leq\displaystyle\lim_{n}\inf\frac{E_n\left((v_n, (q^b_1, q^b_2, 0))\right)}{(h^a_n)^2}=E\left((q^a, q^b)\right).
\end{align}
On the other hand, the properties of $\hat{p}^a, \hat{p}^b$ and the properties of $z^a_1, z^a_2, , z^b_2$ and $z^b_3$ in (\ref{eq:61}) give
\begin{equation}\label{eq:68}
\begin{cases}
\displaystyle\int_{\Omega^a}\left|\frac{\partial z^a_1}{\partial x_1}+\frac{\partial z^a_2}{\partial x_2}+\text{div}\,\hat{p}^a\right|^2dx=|\Theta|\int_{0}^{1}|\text{div}\,\hat{p}^a|^2dx_3+\int_{\Omega^a}\left|\frac{\partial z^a_1}{\partial x_1}+\frac{\partial z^a_2}{\partial x_2}\right|^2dx, &\\
\displaystyle\int_{\Omega^b}\left|\frac{\partial z^b_3}{\partial x_3}+\text{div}\,\hat{p}^b\right|^2dx=\int_{\Theta}|\text{div}\,\hat{p}^b|^2dx_1+\int_{\Omega^b}\left|\frac{\partial z^b_3}{\partial x_3}\right|^2dx. &
\end{cases}
\end{equation} 
Then, 
\begin{align}
&\nonumber\min\{1,\beta, \ell, \beta\ell\}\left[\displaystyle\int_{\Omega^a}\left(\left|\frac{\partial z^a_3}{\partial x_2}\right|^2+\left|\frac{\partial z^a_3}{\partial x_1}\right|^2+\left|\frac{\partial z^a_2}{\partial x_1}-\frac{\partial z^a_1}{\partial x_2}\right|^2+\left|\frac{\partial z^a_1}{\partial x_1}+\frac{\partial z^a_2}{\partial x_2}\right|^2\right)\right.dx\\
&\nonumber\quad+\displaystyle\int_{\Omega^b}\left.\left(\left|\frac{\partial z^b_2}{\partial x_3}\right|^2+\left|\frac{\partial z^b_1}{\partial x_3}\right|^2+\left|\frac{\partial z^b_3}{\partial x_3}\right|^2\right)dx\right]+E((\hat{p}^a, \hat{p}^b))\\
&\label{eq:69}\leq \lim_{n}\inf\frac{E_n\left((p^a_n, p^b_n)\right)}{(h^a_n)^2}\leq \lim_{n}\sup\frac{E_n\left((p^a_n, p^b_n)\right)}{(h^a_n)^2}\leq E\left((q^a, q^b)\right),\quad\text{for all}\, (q^a, q^b)\in P,
\end{align}
from which, thanks to properties of $z^a$ and $z^b$,
\begin{equation}\label{eq:70}
\begin{cases}
\displaystyle \frac{\partial z^a_2}{\partial x_1}-\frac{\partial z^a_1}{\partial x_2}=0 &\text{a.e. in}\quad \Omega^a,\\
\displaystyle \frac{\partial z^a_1}{\partial x_1}+\frac{\partial z^a_2}{\partial x_2}=0 &\text{a.e. in}\quad \Omega^a,\\
\displaystyle \frac{\partial z^a_3}{\partial x_2}=\frac{\partial z^a_3}{\partial x_1}=0 &\text{a.e. in}\quad \Omega^a,\\
\displaystyle \frac{\partial z^b_2}{\partial x_3}=0 &\text{a.e. in}\quad \Omega^b,\\
\displaystyle \frac{\partial z^b_1}{\partial x_3}=0 &\text{a.e. in}\quad \Omega^b,\\
\displaystyle \frac{\partial z^b_3}{\partial x_3}=0 &\text{a.e. in}\quad \Omega^b.
\end{cases}
\end{equation}
In particular $z^a_3=0$ a.e. in $\Omega^a$ and $z^b=0$ a.e. in $\Omega^b$ since $z^a_3\in L^2\left(H^1_m(\Theta),]0, 1[\right)$ and $z^b\in (L^2\left(\Theta, H^1_{m}\left(]-1, 0[\right)\right))^2\times L^2\left(\Theta, H^1_{m}\left(]-1, 0[\right)\right)$. Consequently, inserting (\ref{eq:70}) in (\ref{eq:69}), one obtains that $(\hat{p}^a, \hat{p}^b)$ solves (\ref{eq:23}), and convergence (\ref{eq:24}) holds true. Note that convergences in (\ref{eq:24}) and (\ref{eq:62}) holds true for the whole sequence since the limits are uniquely identified. Moreover, (\ref{eq:22}) follows from (\ref{eq:29}) and Proposition \ref{prop:3}. Now, it remains to prove that convergences in (\ref{eq:29}) and (\ref{eq:62}) are strong. To this end, (\ref{eq:24}), (\ref{eq:14}), (\ref{eq:22}) and (\ref{eq:29}) imply
\begin{align*}
&\displaystyle\lim_{n}\left(\int_{\Omega^a}\left(\beta|\text{rot}\,_n^a p^a_n|^2+|\text{div}\,_n^a p^a_n|^2\right)dx+\frac{h^b_n}{(h^a_n)^2} \int_{\Omega^b}\left(\beta|\text{rot}\,_n^b p^b_n|^2+|\text{div}\,_n^b p^b_n|^2\right)dx\right)\\
&=\displaystyle\left(\int_{\Omega^a}\left(\beta|\text{rot}\,(0, 0, \hat{p}^a)|^2+|\text{div}\,(0, 0, \hat{p}^a)|^2\right)dx+\ell\int_{\Omega^b}\left(\beta|\text{rot}\,(\hat{p}^b_1, \hat{p}^b_2, 0)|^2+|\text{div}\,(\hat{p}^b_1, \hat{p}^b_2, 0)|^2\right)dx\right),
\end{align*}
from which, using (\ref{eq:29}), (\ref{eq:62}) and (\ref{eq:70}), one deduce
\begin{equation*}
\begin{cases}
\left(\text{rot}\,^a_n p^a_n, \text{rot}\,^b_n p^b_n \right)\rightarrow \left(\text{rot}\,(0, 0, \hat{p}^a), \text{rot}\,(\hat{p}^b_1, \hat{p}^b_2, 0) \right) & \text{strongly in}\, \left(L^2(\Omega^a)\right)^3\times\left(L^2(\Omega^b)\right)^3,\\
\left(\text{div}\,^a_n p^a_n, \text{div}\,^b_n p^b_n \right)\rightarrow \left(\text{div}\,(0, 0, \hat{p}^a), \text{div}\,(\hat{p}^b_1, \hat{p}^b_2, 0) \right) & \text{strongly in}\, \left(L^2(\Omega^a)\right)\times\left(L^2(\Omega^b)\right).
\end{cases}
\end{equation*}
Consequently, recalling (\ref{eq:15}), one has
\begin{align}
||D^a_n p^a_n||^2_{\left(L^2(\Omega^a)\right)^9}+\frac{h^b_n}{(h^a_n)^2}||D^b_n p^b_n||^2_{\left(L^2(\Omega^b)\right)^9}&\nonumber=||\text{rot}\,_n^a p^a_n||^2_{\left(L^2(\Omega^a)\right)^3}+||\text{div}\,^a_n p^a_n||^2_{L^2(\Omega^a)}\\
&\label{eq:71}\quad+\frac{h^b_n}{(h^a_n)^2}||\text{rot}\,_n^b p^b_n||^2_{\left(L^2(\Omega^b)\right)^3}+\frac{h^b_n}{(h^a_n)^2}||\text{div}\,^b_n p^b_n||^2_{L^2(\Omega^b)},
\end{align}
and so, as $n\to +\infty$,
\begin{align*}
||D^a_n p^a_n||^2_{\left(L^2(\Omega^a)\right)^9}+\frac{h^b_n}{(h^a_n)^2}||D^b_n p^b_n||^2_{\left(L^2(\Omega^b)\right)^9}&\nonumber\rightarrow||\text{rot}\,(0, 0, \hat{p}^a)||^2_{\left(L^2(\Omega^a)\right)^3}+||\text{div}\,(0, 0, \hat{p}^a)||^2_{L^2(\Omega^a)}\\
&\quad+\ell||\text{rot}\,(\hat{p}^b_1, \hat{p}^b_2, 0)||^2_{\left(L^2(\Omega^b)\right)^3}+\ell||\text{div}\,(\hat{p}^b_1,\hat{p}^b_2 , 0)||^2_{L^2(\Omega^b)}\\
&=||D(0, 0, \hat{p}^a)||^2_{\left(L^2(\Omega^b)\right)^9}+\ell||D(\hat{p}^b_1, \hat{p}^b_2, 0)||^2_{\left(L^2(\Omega^b)\right)^9},
\end{align*}
where the boundary conditions of $(\hat{p}^a, \hat{p}^{b})$ play a crucial role in the last inequatilty. Finally, combining (\ref{eq:71}) with (\ref{eq:29}), (\ref{eq:62}) and (\ref{eq:70}), one obtains that the convergences (\ref{eq:29}) and (\ref{eq:62}) are strong.
\lqqd

\vspace{1.0 cm}

\section{The proof in the case $\ell=0$}\label{cha:5}
\vspace{1.0 cm}

\subsection{A priori estimates on polarization}
\vspace{0.5 cm}
The same arguments used to proving  Proposition \ref{prop:1} give the following result.

\begin{prop}\label{prop:1*}
Assume $(\ref{eq:1})$ with $\ell=0$,  and (\ref{eq:14}). For every $n\in\NN$, let $(p^a_n, p^b_n)$ be a solution of (\ref{eq:13}). Then, there exist a constant $c$ such that
\begin{align}
&\label{eq:25*}||p^a_{n}||_{\left(L^4(\Omega^a)\right)^3}\leq c,\quad \left|\left|\frac{\sqrt[4]{h^b_n}}{\sqrt{h^a_n}}p^b_n\right|\right|_{\left(L^4(\Omega^b)\right)^3}\leq c\quad \text{for all}\,\, n\in\NN,\\ 
 &\label{eq:26*}||D^a_np^a_{n}||_{\left(L^2(\Omega^a)\right)^9}\leq c,\quad \left|\left|\frac{\sqrt{h^b_n}}{h^a_n}D^b_np^b_n\right|\right|_{\left(L^2(\Omega^b)\right)^9}\leq c\quad \text{for all}\,\, n\in\NN.               
 \end{align}
 
 \vspace{0.5 cm}
 
 \dem
 Function $0$ belonging to $P_n$ gives
\begin{align}
&\nonumber\displaystyle\int_{\Omega^a}\left(\beta|\text{rot}_n^a\, p^a_n|^2+|\text{div}_n^a\, p^a_n|^2+\alpha\left(|p^a_n|^4-2|p^a_n|^2\right)+|D_n^a\phi^a_{(p^a_n, p^b_n)}|^2\right)dx\\
&\nonumber\quad +\frac{h^b_n}{(h^a_n)^2}\int_{\Omega^b}\left(\beta|\text{rot}_n^b\, p^b_2|^2+|\text{div}_n^b\, p^b_n|^2+\alpha\left(|p^b_n|^4-2|p^b_n|^2\right)^2+|D_n^b\phi^b_{(p^a, p^b)}|^2\right)dx,\\
&\label{eq:27*}\quad\leq\displaystyle\frac{1}{2}\int_{\Omega^a}\left(|f^a_n|^2+|p^a_n|^2\right)dx+\frac{h^b_n}{(h^a_n)^2}\displaystyle\frac{1}{2}\int_{\Omega^b}\left(|f^b_n|^2+|p^b_n|^2\right)dx\quad \text{for all $n\in\NN$}.
\end{align}
Estimates (\ref{eq:27*}) implies
\begin{align*}
&\nonumber\displaystyle\int_{\Omega^a}\alpha\left(|p^a_n|^4-\left(2+\frac{1}{2\alpha}\right)|p^a_n|^2\right)dx+\frac{h^b_n}{(h^a_n)^2}\int_{\Omega^b}\alpha\left(|p^b_n|^4-\left(2+\frac{1}{2\alpha}\right)|p^b_n|^2\right)dx\\
&\quad\leq\displaystyle\frac{1}{2}\int_{\Omega^a}\left(|f^a_n|^2+|p^a_n|^2\right)dx+\frac{h^b_n}{(h^a_n)^2}\displaystyle\frac{1}{2}\int_{\Omega^b}\left(|f^b_n|^2+|p^b_n|^2\right)dx\quad \text{for all $n\in\NN$},
\end{align*}
which gives
\begin{align}
&\nonumber\displaystyle\int_{\Omega^a}\alpha\left(|p^a_n|^2-\left(1+\frac{1}{4\alpha}\right)\right)^2dx+\frac{h^b_n}{(h^a_n)^2}\int_{\Omega^b}\alpha\left(|p^b_n|^2-\left(1+\frac{1}{4\alpha}\right)\right)^2dx\\
&\label{eq:28*}\quad\leq\displaystyle\alpha\left(1+\frac{1}{4\alpha}\right)\left(|\Omega^a|+\frac{h^b_n}{(h^a_n)^2}|\Omega^b|\right)+\frac{1}{2}\int_{\Omega^a}|f^a_n|^2dx+\frac{h^b_n}{(h^a_n)^2}\displaystyle\frac{1}{2}\int_{\Omega^b}|f^b_n|^2dx\quad \text{for all $n\in\NN$}.
\end{align}
Then the estimate $(\ref{eq:25*})$ follow from (\ref{eq:28*}), (\ref{eq:1}) with $\ell=0$ and (\ref{eq:14}). The estimate (\ref{eq:26*}) follow from (\ref{eq:27*}) with $\ell=0$, (\ref{eq:14}), (\ref{eq:25*}), the continuous embedding of $L^4$ into $L^2$ and (\ref{eq:15}).
\lqqd
\end{prop}

\vspace{1.0 cm}

\begin{coro}\label{coro:1*} 
Assume (\ref{eq:1}) with $\ell=0$ and (\ref{eq:14}). For every $n\in\NN$, let $(p^a_n, p^b_n)$ be a solution of (\ref{eq:13}), let $P_0$ be defined in (\ref{eq:16*}). Then 
\begin{equation}\label{eq:29*}
          \frac{\sqrt{h^b_n}}{h^a_n}p^b_{n}\rightharpoonup 0 \, \text{weakly in} \, \left(H^1(\Omega^b)\right)^3\,  \text{ and strongly in}\,  \left(L^4(\Omega^b)\right)^3.
\end{equation}
Moreover, there exist a subsequence of $\NN$, still denoted by $\{n\}$, and (in possible dependence on the subset) $\hat{p}^a\in P_0$ such that
\begin{equation}\label{eq:9999*}
    p^a_{n}\rightharpoonup (0,0,\hat{p}^a) \, \text{weakly in} \, \left(H^1(\Omega^a)\right)^3\,  \text{ and strongly in}\,  \left(L^4(\Omega^a)\right)^3.
\end{equation}

\end{coro}

\vspace{0.5 cm}

\dem Proposition \ref{prop:1*} implies (\ref{eq:29*}) and the existence of  a subsequence of $\NN$, still denoted by $\{n\}$, and (in possible dependence on the subsequence) $(\hat{p}^a_1, \hat{p}^a_2, \hat{p}^a)\in \left(H^1(\Omega^a)\right)^3$ independent of $x_1$ and $x_2$,  such that
\begin{equation}\label{eq:30*}
          p^a_{n}\rightharpoonup (\hat{p}^a_1, \hat{p}^a_2, \hat{p}^a) \, \text{weakly in} \, \left(H^1(\Omega^a)\right)^3\,  \text{ and strongly in}\,  \left(L^4(\Omega^a)\right)^3,       
\end{equation}
and $(\hat{p}^a_1, \hat{p}^a_2, \hat{p}^a)\cdot \nu^a=0$ on $\partial\Omega^a\setminus\left(\Theta\times\{0\}\right)$.

\vspace{0.5 cm}

In particular, this implies
\begin{equation*}
          \hat{p}^a_1=\hat{p}^a_2=0 \, \text{in} \, \Omega^a\,  \text{and}\,\,  \hat{p}^a(1)=0,
\end{equation*}
It remains to prove that
\begin{equation}\label{eq:31*}
 \hat{p}^a(0)=0.
 \end{equation}
 
\vspace{0.5 cm}

The trace of $p^b_{n, 3}$ vanishing on $\Theta\times\{-1\}$ implies
\begin{equation*}
\displaystyle\left|\frac{1}{\sqrt{h^b_n}h^a_n}p^b_{n, 3}(x', 0)\right|^2\leq \int_{-1}^{0}\left|\frac{1}{\sqrt{h^b_n}h^a_n}p^b_{n, 3}(x', t)\right|^2dt\quad\text{in}\,\Theta\quad\text{for all}\,n\in\NN.
\end{equation*}
Integrating this inequality over $\Theta$ and using the second estimate in (\ref{eq:26*}) gives
\begin{equation}\label{eq:32*}
\frac{1}{\sqrt{h^b_n}h^a_n}p^b_{n, 3}(\cdot, \cdot, 0)\to 0\quad\text{strongly in}\, L^2(\Theta).
\end{equation}
Finally, (\ref{eq:31*}) follows from (\ref{eq:30*}), the junction condition in (\ref{eq:8}), and (\ref{eq:32*}). Indeed,
\begin{align*}
\displaystyle|\Theta||\hat{p}^a(0)|^2&=\lim_{n}\int_{\Theta}|p^a_{n, 3}(x', 0)|^2dx'\\
&=\lim_{n}\int_{\Theta}|p_{n, 3}^b(h^a_n x', 0)|^2dx'\\
&=\lim_{n}\int_{\Theta}\left|\frac{1}{h^a_n}p_{n, 3}^b(x', 0)\right|^2dx'\\
&=\lim_{n}h^b_n\cdot \lim_{n}\int_{\Theta}\left|\frac{1}{\sqrt{h^b_n}h^a_n}p_{n, 3}^b(x', 0)\right|^2dx'=0.
\end{align*}
 \lqqd

\vspace{0.5 cm}
\subsection{A convergence result for problem (\ref{eq:11})}
\vspace{0.5 cm}

\begin{prop}\label{prop:3*}
Assume (\ref{eq:1}) with $\ell=0$. Let $\{(q^a_n, q^b_n)\}_{n\in\NN}\subset \left(L^2(\Omega^a)\right)^3\times (L^2(\Omega^b))^3$, and let $q^a\in L^2(\Omega^a)$ such that $q^a$ is independent of $(x_1,x_2)$, and 
\begin{equation}\label{eq:45*}
\left(q^a_n, \frac{\sqrt{h^b_n}}{h^a_n}q^b_n\right)\rightarrow ((0, 0, q^a), 0)\quad\text{strongly in}\, \left(L^2(\Omega^a)\right)^3\times (L^2(\Omega^b))^3. 
\end{equation}
Moreover, for $n\in\NN$ let $(\phi^a_{(q^a_n, q^b_n)}, \phi^b_{(q^a_n, q^b_n)} )$ be the unique solution of 
\begin{equation}\label{eq:46*}
\begin{cases} 
      \left(\phi^a_{(q^a_n, q^b_n)},\phi^b_{(q^a_n, q^b_n)} \right)\in U_n, \quad \displaystyle\int_{\Omega^a}\phi^a_{(q^a_n, q^b_n)}dx=0, &   \\
      \displaystyle\int_{\Omega^a}\left(\left(-D^a_n\phi^a_{(q^a_n, q^b_n)}+q^a_n\right)\cdot D^a_n \phi^a\right)dx+\frac{h^b_n}{(h^a_n)^2}\displaystyle\int_{\Omega^b}\left(\left(-D^b_n\phi^b_{(q^a_n, q^b_n)}+q^b_n\right)\cdot D^b_n \phi^b\right)dx=0, & 
         \end{cases}
\end{equation}
$(\phi^a, \phi^b)\in U_n,$ where $U_n$ is defined in (\ref{eq:9}). Furthermore, let $\psi^a_{q^a}$ the unique solution of (\ref{eq:19*}). Then one has
\begin{equation}\label{eq:47*}
\begin{cases} 
      \displaystyle\left(\phi^a_{(q^a_{n},q^b_{n})}, \frac{\sqrt{h^b_n}}{h^a_n}\phi^b_{(q^a_{n}, q^b_{n})} \right)\rightarrow(\psi^a_{q^a}, 0) & \text{strongly in}\, \left(H^1(\Omega^a)\right)\times \left(H^1(\Omega^b)\right), \\
     \displaystyle\left(\frac{1}{h^a_n}\frac{\partial \phi^a_{(q^a_{n}, q^b_{n})}}{\partial x_1}, \frac{1}{h^a_n}\frac{\partial \phi^a_{(q^a_{n}, q^b_{n})}}{\partial x_2} \right)\rightarrow(0,0) & \text{strongly in}\, \left(L^2(\Omega^a)\right)^2,\\
       \displaystyle\frac{1}{h^b_n}\frac{\partial  \left( \frac{\sqrt{h^b_n}}{h^a_n}\phi^b_{(q^a_{n}, q^b_{n})}\right) }{\partial x_3}\rightarrow 0 & \text{strongly in}\,  L^2(\Omega^b).  
 \end{cases}
\end{equation}
\end{prop}

\vspace{0.5 cm}

\dem
In this proof, $c$ denotes any positive constant independent of $n\in\NN$. Choosing $(\phi^a, \phi^b)=\left(\phi^a_{(q^a_{n}, q^b_{n})}, \phi^b_{(q^a_{n}, q^b_{n})} \right)$ as test function in (\ref{eq:46*}) applying Young inequality, and using (\ref{eq:45*}) give
\begin{equation}\label{eq:48*}
||D^a_n\phi^a_{(q^a_{n}, q^b_{n})}||_{(L^2(\Omega^a))^3}\leq c,\quad \left|\left|\frac{\sqrt{h^b_n}}{h^a_n}D^b_n\phi^b_{(q^a_{n}, q^b_{n})}\right|\right|_{(L^2(\Omega^b))^3}\leq c\quad \text{for all}\, n\in\NN.
\end{equation}
the first estimate in (\ref{eq:48*}) implies
\begin{equation}\label{eq:49*}
 ||\phi^a_{(q^a_{n}, q^b_{n})}||_{H^1(\Omega^a)}\leq c\quad \text{for all}\, n\in\NN,
\end{equation}
since $\displaystyle\int_{\Omega^a}\phi^a_{(q^a_{n}, q^b_{n})}dx=0$. The next step is devoted to proving
\begin{equation}\label{eq:50*}
 \left|\left|\frac{\sqrt{h^b_n}}{h^a_n}\phi^b_{(q^a_{n}, q^b_{n})}\right|\right|_{H^1(\Omega^b)}\leq c\quad \text{for all}\, n\in\NN.
\end{equation}
The junction condition in (\ref{eq:9}) gives  
\begin{align}
\displaystyle\int_{h^a_n\Theta}\left|\frac{\sqrt{h^b_n}}{h^a_n}\phi^b_{(q^a_{n}, q^b_{n})}(x', 0)\right|^2 dx' &\nonumber=\frac{h^b_n}{(h^a_n)^2}(h^a_n)^2\int_{\Theta}\left|\phi^b_{(q^a_{n}, q^b_{n})}( h^a_n x', 0)\right|^2 dx'\\
&\label{eq:51*}=h^b_n\int_{\Theta}\left|\phi^a_{(q^a_{n}, q^b_{n})}( x', 0)\right|^2 dx'
\end{align}
$\,\text{for all}\, n\in\NN$.
Then, (\ref{eq:51*}) and (\ref{eq:49*}) provide
\begin{equation*}
\left|\left|\frac{\sqrt{h^b_n}}{h^a_n}\phi^b_{(q^a_{n}, q^b_{n})}\right|\right|_{L^2\left(h^a_n\Theta\times\{0\}\right)}\leq \sqrt{h^b_n}c,\quad\text{for all}\,n\in\NN,
\end{equation*}
which implies
\begin{equation}\label{eq:52*}
\left|\left|\frac{\sqrt{h^b_n}}{h^a_n}\phi^b_{(q^a_{n}, q^b_{n})}\right|\right|_{L^2\left(h^a_n\Theta\times]-1, 0[\right)}\leq \sqrt{h^b_n}c,\quad\text{for all}\,n\in\NN,
\end{equation}
by virtue of the (\ref{eq:48*}) ensures
\begin{equation*}
\left|\left|\frac{\sqrt{h^b_n}}{h^a_n}\phi^b_{q_n}\right|\right|_{L^2\left(\{0\}\times\{0\}\times]-1,0[\right)}\leq c,\quad\text{for all}\,n\in\NN,
\end{equation*}
which combined again with the second estimates in (\ref{eq:48*}) proves (\ref{eq:50*}). The estimates (\ref{eq:48*}), (\ref{eq:49*}), and (\ref{eq:50*}) ensure the existence of a subsequence of $\NN$, still denotes by $\{n\}$ and (in possible dependence on the subsequence) $(\tau^a, \tau^b)\in H^1(\Omega^a)\times H^1(\Omega^b)$ with $\tau^a$ independent of $(x_1, x_2)$ and $\tau^b$ independent of $x_3$  , $(\xi^a, \zeta^a)\in (L^2(\Omega^a))^2$ and $\zeta^b\in L^2(\Omega^b)$ such that
\begin{align} 
      \displaystyle\left(\phi^a_{(q^a_{n},q^b_{n})}, \frac{\sqrt{h^b_n}}{h^a_n}\phi^b_{(q^a_{n}, q^b_{n})} \right)\rightharpoonup(\tau^a, \tau^b) &\label{eq:53*}\, \text{weakly in}\, \left(H^1(\Omega^a)\right)\times \left(H^1(\Omega^b)\right), \\
     \displaystyle\left(\frac{1}{h^a_n}\frac{\partial \phi^a_{(q^a_{n}, q^b_{n})}}{\partial x_1}, \frac{1}{h^a_n}\frac{\partial \phi^a_{(q^a_{n}, q^b_{n})}}{\partial x_2} \right)\rightharpoonup (\xi^a,\zeta^a) &\label{eq:54*}\, \text{weakly in}\, \left(L^2(\Omega^a)\right)^2,\\
       \displaystyle  \frac{1}{h^b_n}\frac{\partial \left(\frac{\sqrt{h^b_n}}{h^a_n}  \phi^b_{(q^a_{n}, q^b_{n})} \right) }{\partial x_3}\rightharpoonup \zeta^b &\label{eq:55*}\, \text{weakly in}\,  L^2(\Omega^b),  
 \end{align}
and
\begin{equation}\label{eq:56*}
\displaystyle\int_{0}^{1}\tau^a dx_3=0.
\end{equation}
The next step is devoted to proving 
\begin{equation}\label{eq:57*}
\tau^b(0')=0.
\end{equation}

\vspace{0.5 cm}

Indeed, the junction condition in (\ref{eq:9}) gives
\begin{align}
\displaystyle\int_{\Theta}&\nonumber\frac{\sqrt{h^b_n}}{h^a_n}\phi^a_{(q^a_n, q^b_n)}(x', 0)\varphi(x')dx'\\
&\label{eq:58*}=\displaystyle\int_{\Theta}\frac{\sqrt{h^b_n}}{h^a_n}\phi^b_{(q^a_n, q^b_n)}(h^ax', 0)\varphi(x')dx'\quad\text{for all}\,\,\varphi\in C^{\infty}_0(\Theta).
\end{align}
Moreover, (\ref{eq:1}) with $\ell=0$ and the convergence of the first term in (\ref{eq:53*}) imply
\begin{equation}\label{eq:59*}
\frac{\sqrt{h^b_n}}{h^a_n}\phi^a_{(q^a_n, q^b_n)}(\cdot, \cdot, 0)\to 0\quad\text{strongly in}\, L^2(\Theta),
\end{equation}
while, using the second convergence in (\ref{eq:53*}) and (\ref{eq:54*}), again (\ref{eq:1}) with $\ell=0$, and arguing as in \cite{REF20} proof of Proposition 5.4, one can prove
\begin{align}
\displaystyle\int_{\Theta}&\nonumber\frac{\sqrt{h^b_n}}{h^a_n}\phi^b_{(q^a_n, q^b_n)}(h^ax', 0)\varphi(x')dx'\\
&\label{eq:60*}=\displaystyle\int_{\Theta}\tau^b(x',0)\varphi(x')dx'\quad\text{for all}\,\,\varphi\in C^{\infty}_0(\Theta).
\end{align}
Then (\ref{eq:57*}) follows from (\ref{eq:58*}), (\ref{eq:59}), and (\ref{eq:60*}).

\vspace{0.5 cm}

To identify $\tau^a=\psi^a_{q^a}$, it is enough to pass to the limits, as $n$ diverges, in the equation in (\ref{eq:46*}) with a test functions $(\phi^a, \phi^b)$ such that
\begin{align*}
\phi^a(x)&=\displaystyle\frac{1}{h^a_n}\psi^a(x_3)\,\,\text{if}\,\,x=(x_1, x_2, x_3)\in\Omega^a,\\
\phi^b(x)&=\displaystyle\frac{1}{h^a_n}\psi^a(0)\,\,\text{if}\,\,x=(x_1, x_2, x_3)\in\Omega^b,
\end{align*}
with $\psi^a\in C^{\infty}([0,1])$, to use (\ref{eq:1}) with $\ell=0$, previous convergences, the density of the space $C^{\infty}([0,1])$ in $H^1([0,1])$, and (\ref{eq:56*}).
\vspace{0.5 cm}

To identify $\tau^b=0$, it is enough to pass to the limit, as $n$ diverges, in the equation (\ref{eq:46*}) with test functions $(\phi^a, \phi^b)$ such that
\begin{align*}
\phi^a(x)&=0\,\,\text{if}\,\,x=(x_1, x_2, x_3)\in\Omega^a,\\
\phi^b(x)&=\displaystyle\frac{1}{\sqrt{h^b_n}h^a_n}\psi^b(x')\,\,\text{if}\,\,x=(x', x_3)\in\Omega^b,
\end{align*}
with
\begin{equation*}
\psi^b\in A=\left\{\nu\in C^{\infty}(\overline{\Theta}): \nu=0\,\,\text{in}\,\, \Theta\right\},
\end{equation*}
the previous convergences, the density of $A$ in the space $H^1(\Theta)$ with zero trace, and to take into account (\ref{eq:57*}).
To identifty $\xi^a, \zeta^a$ and $\zeta^b$ and to prove that all the previous convergences are strong, one can argue as in the last part of the proof of Proposition \ref{prop:3}.
\lqqd

\vspace{1.0 cm}
\subsection{The proof of Theorem \ref{teo:2}}
\vspace{0.5 cm}

\subsubsection*{Proof of Theorem \ref{teo:2}} Corollary \ref{coro:1*} asserts that (\ref{eq:30*}) holds true and that there exist a subsequence of $\NN$, still denoted by $\{n\}$, and (in possible dependence on the subsequence) of $\hat{p}^a\in P_0$ satisfying (\ref{eq:31*}). Moreover, once can prove the existence of a subsequence of $\NN$, still denoted by $\{n\}$, and  (in possible dependence on the subsequence) of
\begin{equation*}
z^a\in L^2(H^1_0(\Theta), ]0,1[)\times (L^2(H^1_m(\Theta), ]0,1[))^2\\
\end{equation*}
\ such that
\begin{equation}\label{eq:62*}
\begin{cases} 
      &  \displaystyle\frac{1}{h^a_n}D_{x'} p^a_n \rightharpoonup D_{x'} z^a\,\, \text{weakly in}\, (L^2(\Omega^a))^6,\\
       &  \displaystyle\frac{1}{h^b_n}\frac{\partial\left(\frac{\sqrt{h^b_n}}{h^a_n} p^b_n\right)}{\partial x_3} \rightharpoonup \zeta\,\, \text{weakly in}\, (L^2(\Omega^b))^3.
        \end{cases}
\end{equation}
The next step is to identify $\hat{p}^a$, $z^a$ and $\zeta$. Let
\begin{equation}\label{eq:16*}
P_0^{\text{reg}}=\left\{q^a\in H^1([0, 1]):\,\,\text{for some}\,\,\delta>0\,\,(\text{depending on}\,q^a),\,\,q^a=0\,\,\text{in}\,\, [0, \delta]\cup[1-\delta, 1] \right\}.
\end{equation}

Indeed, the boundary condition on $p^a_n$ and the Poincar\'e inequality give, for $i=1,2$,
\begin{equation}\label{eq:63*}
\displaystyle\left|\left|\frac{1}{h^a_n}p^a_{n,i}(\cdot, \cdot, x_3)\right|\right|_{H^1(\Theta)}\leq \frac{c}{h^a_n}\left|\left|\left( \frac{\partial p^a_{n,i}(\cdot, \cdot, x_3) }{\partial x_1}, \frac{\partial p^a_{n,i}(\cdot, \cdot, x_3)}{\partial x_2}\right)\right|\right|_{L^2(\Theta)}
\end{equation}
for a.e. $x_3\in]0,1[$ and all $n\in\NN$, where $c$ is the Poincar\'e constant in $H^1_0(\Theta)$, while, for $i=3$, the Poincar\'e- Wirtinger inequality gives
\begin{align}
&\nonumber\displaystyle\left|\left|\frac{1}{h^a_n}\left(p^a_{n,3}(\cdot, \cdot, x_3)-\int_{\Theta}p^a_{n, 3}(x_1, x_2, x_3)dx_1dx_2\right)\right|\right|_{H^1_m(\Theta)}\\
&\label{eq:64*}\displaystyle\quad\leq \frac{c'}{h^a_n}\left|\left|\left( \frac{\partial p^a_{n,3}(\cdot, \cdot, x_3) }{\partial x_1}, \frac{\partial p^a_{n,3}(\cdot, \cdot, x_3)}{\partial x_2}\right)\right|\right|_{L^2(\Theta)}
\end{align}
for a.e. $x_3\in]0,1[$, and all $n\in\NN$, where the subscript ``$m$'' means zero average, and $c'$ is the Poincar\'e-Wirtinger constant in $H^1_{m}(\Theta)$. Integrating (\ref{eq:63*}) and (\ref{eq:64*}) over $x_3\in]0, 1[$ and using the first estimate in (\ref{eq:26}) imply the first convergence in (\ref{eq:62*}).Similarly, one proves the remaining convergences in $(\ref{eq:62})$. 

The next step is devoted to identifying $\hat{p}^a, \hat{p}^b, z^a$, and $z^b$. To this end, for every
\begin{equation*}
(q^a, q^b)\in P^{\text{reg}}_0,
\end{equation*}
where $P^{\text{reg}}_0$ is defined in (\ref{eq:16*}), consider a sequence $\{v_n\}_{n\in\NN}\subset (H^1(\Omega^a))^3$ (depending on $(q^a, q^b)$) such that
\begin{equation}\label{eq:65*}
\begin{cases} 
      \displaystyle\left( v_n, (q^b, 0, 0) \right)\in P_n & \text{for all}\, n\in\NN, \\
     \displaystyle v_n \rightarrow (0, 0, q^a) & \text{strongly in}\, (L^4(\Omega^a))^3,\\
       \displaystyle\left(\frac{1}{h^a_n}D_{x_1} v_n, \frac{1}{h^a_n}D_{x_2} v_n, D_{x_3} v_n\right)\rightarrow(0,0, D_{x_3} q^a) & \text{strongly in}\, \left(L^2(\Omega^a)\right)^9.  
 \end{cases}
\end{equation}
For instance, setting
\begin{align*}
v_n(x)&=(0, 0, q^a)(x_3)\quad\text{if}\, x=(x_1, x_2, x_3)\in \Theta\times]h^a_n, 1[\,\text{and}\\
v_n(x)&=\left(q^b(h^a_nx_1)\frac{h^a_n-x_3}{h^a_n}, 0, q^a(x_3)\right)\quad\text{if}\, x=(x_1, x_2, x_3)\in \Theta\times[0, h^a_n],
\end{align*}
the first two properties in (\ref{eq:65*}) can be immediately verified by virtue of the properties of $P_{\text{reg}}$, while the last ones follows from
\begin{align*}
\displaystyle\int_{\Theta\times]0, h_n[} \left|\frac{1}{h^a_n}D_{x_1}v_n\right|^2dx&=\int_{\Theta\times]0, h^a_n[}\left|D_{x_1} q^b(h^a_n x_1)\left(\frac{h^a_n-x_3}{h^a_n}\right)\right|^2dx\\
&\displaystyle=\int_{\Theta}\left|D_{x_1} q^b(h^a_n x_1)\right|^2d{x_1}dx_{2}\int_{]0, h_n[}\left(1-\frac{x_3}{h^a_n}\right)^2dx_3\\
&\leq ||q^b||^2_{W^{1,\infty}(]-\frac{1}{2},\frac{1}{2}[^2)}h^a_n,\quad\text{for all}\, n\in\NN.
\end{align*}
\begin{equation*}
\displaystyle\int_{\Theta\times]0, h^a_n[} \left|\frac{1}{h^a_n}D_{x_2}v_n\right|^2dx=0,\quad\text{for all}\, n\in\NN.
\end{equation*}
\begin{align*}
\displaystyle\int_{\Theta\times]0, h^a_n[} \left| D_{x_3} v_n\right|^2dx&=\int_{\Theta\times]0, h^a_n[}\left|q^b(h^a_nx_1)\frac{1}{h^a_n}\right|^2dx+\int_{\Theta\times]0, h^a_n[} \left|D_{x_3}q^{a}(x_3)\right|^2dx\\
&\leq\int_{\Theta\times]0, h^a_n[}\left|\frac{1}{h^a_n}\left(q^b(h^a_nx_1)-q^b(0)\right)\right|^2dx\\
&\quad +\int_{\Theta\times]0, h^a_n[} \left|D_{x_3}q^{a}(x_3)\right|^2dx\\
&\leq\left(||q^b||^2_{W^{1,\infty}(]-\frac{1}{2},\frac{1}{2}[)}+||q^a||^2_{W^{1,\infty}(]0, 1[)}\right)h^a_n,\quad\text{for all}\,n\in\NN,
\end{align*}
where again the properties of $P^{\text{reg}}_0$ played a crucial role.

\vspace{0.5 cm}

Now, fixing $(q^a, q^b)\in P^{\text{reg}}_0$ and choosing $(v_n, (q^b, 0, 0))$ as test function in (\ref{eq:13}) with $v_n$ satisfying (\ref{eq:65*}), give
\begin{equation}\label{eq:66*}
\displaystyle\frac{1}{(h^a_n)^2}E_n\left((p^a_n, p^b_n)\right)\leq \frac{1}{(h^a_n)^2}E_n\left(v_n, (q^b, 0, 0)\right)\quad\text{for all}\,n\in\NN.
\end{equation}
Then, passing to the limit in (\ref{eq:66*}), as $n$ diverges, and using (\ref{eq:7}), (\ref{eq:29*}), (\ref{eq:62*}), (\ref{eq:65*}), Proposition \ref{prop:3} and an l.s.c. argument imply

\begin{align}
&\nonumber\displaystyle\int_{\Omega^a}\left(\beta\left(\left|\frac{\partial z^a_3}{\partial x_2}\right|^2+\left|\frac{\partial z^a_3}{\partial x_1}\right|^2+\left|\frac{\partial z^a_2}{\partial x_1}-\frac{\partial z^a_1}{\partial x_2}\right|^2\right)+\left(\left|\frac{\partial z^a_1}{\partial x_1}+\frac{\partial z^a_2}{\partial x_2}+\text{div}\,\hat{p}^a\right|^2\right)\right)dx\\
&\nonumber\quad+|\Theta|\displaystyle\int_{0}^{1}\left(\alpha(|\hat{p}^a|^2-1)^2+|D_{x_3}\psi^a_{(\hat{p}^a, \hat{p}^b)}|^2\right)dx_3+\int_{0}^{1}\left(\int_{\Theta}f^a_3 dx_1dx_2\cdot \hat{p}^a \right)dx_3\\
&\nonumber\leq\lim_{n}\inf\frac{E_n\left((p^a_n, p^b_n)\right)}{(h^a_n)^2}
\label{eq:67*}\leq\displaystyle\lim_{n}\inf\frac{E_n\left((v_n, (q^b_1, 0, 0))\right)}{(h^a_n)^2}=E_0\left(q^a\right),\quad \text{forall}\, q^a=(0,0,q^a_3)\,\, \text{s.t.}\, q^a_3\in P_0.
\end{align}
These inequalities hold true also for any $q^a\in P^{\text{reg}}_0$, since $P^{\text{reg}}_0$ is dense in $P_0$. One can conclude the proof arguing as in the last part of the proof of Theorem \ref{teo:1}.
\lqqd

\vspace{1.0 cm}
\section{The case for the boundary condition $\textbf{P} // e_3$}\label{chap:6}
\vspace{0.5 cm}
Another very significant choice for an energetic approach consists in explicitly considering the energetic
contribution for the polarization field given by the gradient term. Now, if we consider the following non-convex and non-local energy associated with $\Omega_{n}$:
\begin{equation}\label{eqjuly161}
\mathcal{S}_{n}:\textbf{P}\in\left(H^1(\Omega_n)\right)^3\to \displaystyle\int_{\Omega_n}\left(|D\textbf{P}|^2+\alpha\left(|\textbf{P}|^2-1\right)^2+|D\varphi_{\textbf{P}}|^2+(\textbf{F}_n\cdot\textbf{P})\right)dx,
\end{equation}
where, $\alpha$ is a positive constant, $\textbf{F}_n\in\left(L^2(\Omega_n)\right)^3$, $\nu$ denotes the unit normal on $\partial\Omega_n$ and $\varphi_{\textbf{P}}\in H^1(\Omega_n)$ is the unique solution, up to an additive constant, of (\ref{eq:3}).
Set 
\begin{equation*}
\mathcal{P}^*_n=\{\textbf{P}\in\left(H^1(\Omega_n)\right)^3:\textbf{P}// e_3\quad\text{on}\,\partial\Omega_n\}.
\end{equation*}
The direct method of Calculus of Variations ensures that the following problems
\begin{equation}\label{eqjuly162}
\min \left\{\mathcal{S}_n(\textbf{P}):\textbf{P}\in\mathcal{P}^*_n\right\}
\end{equation}
and
\begin{equation}\label{eqjuly17}
\min \left\{\mathcal{S}_n(\textbf{P}):\textbf{P}\in\mathcal{P}_n\right\},
\end{equation}
admit solution, where $e_3=(0,0,1)$ and $//$ is the symbol of parallelism and $\mathcal{P}_n$ is defined in (\ref{rescaprobjuly19}). 
\vspace{1.0 cm}

\subsection{Rescaling the problem}

\vspace{0.5 cm}

Similar to previous problems, we can use the rescaling given in (\ref{eq:6}) in order to obtain rescaled versions for  the energy $(\ref{eqjuly161})$ and problem (\ref{eqjuly162}) as follows:
Set
\begin{align}
P_n^*=\left\{(p^a, p^b)\in\left(H^1(\Omega^a)\right)^3\times \left(H^1(\Omega^b\right))^3\right. &\nonumber: p^a// e_3\quad \text{on}\, \partial\Omega^a \setminus(\Theta\times\{0\}) \\
&\nonumber \left. p^b//e_3\quad \text{on}\, \partial\Omega^b \setminus(\Theta\times\{0\})\right.\\
&\nonumber \left. p^b_1=p^b_2=0\quad \text{on}\, \left(\Theta\setminus h^a_n\Theta\right)\times\{0\}\right.\\
&\label{eq:8july16} \left. p^a(x', 0)=p^b(h^a_n x', 0)\quad \text{a.e. in}\, \Theta \right\}.
\end{align}
Then $\mathcal{S}_n$ defined in (\ref{eqjuly161}) is rescaled by
\begin{align}
S_{n}:(p^a, p^b)\in P_{n}\to \displaystyle &\nonumber (h^a_n)^2\int_{\Omega^a}\left(|D_n^a\, p^a|^2+\alpha\left(|p^a|^2-1\right)^2-|D_n^a\phi^a_{(p^a, p^b)}|^2+(f^a_n\cdot p^a)\right)dx\\
&\label{eq:10}\quad +h^b_n\int_{\Omega^b}\left(|D_n^b\, p^b|^2+\alpha\left(|p^b|^2-1\right)^2-|D_n^b\phi^b_{(p^a, p^b)}|^2+(f^b_n\cdot p^b)\right)dx,
\end{align}
where $\left(\phi^a_{(p^a, p^b)},\phi^b_{(p^a, p^b)} \right)$ is the unique solution of (\ref{eq:11}).

Note that if $\textbf{P}_n$ solves (\ref{eqjuly162}), then $(p^a_n, p^b_n)$ defined by
\begin{equation*}
p^a_n(x_1, x_2, x_3)=\textbf{P}_n(h^a_n x', x_3)\quad \text{in}\, \Omega^a,\quad p^b_n(x_1, x_2, x_3)=\textbf{P}_n(x', h^b_n x_3)\quad \text{in}\, \Omega^b,
\end{equation*}
solves
\begin{equation}\label{eq:13july16}
\min \{S_n((p^a, p^b)): (p^a, p^b)\in P^*_n\}.
\end{equation}

Analogously, note that if $\textbf{P}_n$ solves (\ref{eqjuly17}), then $(p^a_n, p^b_n)$ defined by
\begin{equation*}
p^a_n(x_1, x_2, x_3)=\textbf{P}_n(h^a_n x', x_3)\quad \text{in}\, \Omega^a,\quad p^b_n(x_1, x_2, x_3)=\textbf{P}_n(x', h^b_n x_3)\quad \text{in}\, \Omega^b,
\end{equation*}
solves
\begin{equation}\label{eq:13july16new}
\min \{S_n((p^a, p^b)): (p^a, p^b)\in P_n\}.
\end{equation}

Set
\begin{equation}\label{newpjuly16}
P^*=\left\{(q^a_3,q^b_3)\in H^1(]0, 1[)\times H^1(\Theta),\, q^a_3(0)=q^{b}_3(0') \right\},
\end{equation}
one can prove the results shown below by closely following the ideas given in \cite{REF24*}:

\begin{teo}\label{teo:1july16}
Assume (\ref{eq:1}) with $\ell\in]0,+\infty[$ and (\ref{eq:14}). For every $n\in \NN$, let $(p^a_n, p^b_n)$ a solution of $(\ref{eq:13july16})$, and $\left(\phi^a_{(p^a_n, p^b_n)}, \phi^b_{(p^a_n, p^b_n)}\right)$ be the unique solution of $(\ref{eq:11})$ with $(p^a, p^b)=(p^a_n, p^b_n)$. Moreover, let $P^*$ and $E$ defined by (\ref{newpjuly16}) and (\ref{eq:18}), (\ref{eq:19}) respectively. Then there exist an increasing sequence of positive integer numbers $\{n_i\}_{i\in\NN}$ and (in possible dependence on the subsequence) $(\hat{p}^a_3,\hat{p}^b_3)\in P^*$ such that
\begin{equation}\label{eq:20july16}
\begin{array}{rcl}
          p^a_{n_i}\rightarrow (0, 0, \hat{p}^a_3) & \text{strongly in} & \left(H^1(\Omega^a)\right)^3\,  \text{ and strongly in}\,  \left(L^4(\Omega^a)\right)^3,\\
          p^b_{n_i}\rightarrow (0, 0, \hat{p}^b_3) & \text{strongly in} & \left(H^1(\Omega^b)\right)^3\,  \text{ and strongly in}\,  \left(L^4(\Omega^b)\right)^3,
\end{array}
\end{equation}
(\ref{eq:21})-(\ref{eq:22}) hold and where $(\hat{p}^a, (\hat{p}^b_1, \hat{p}^b_2))$ solves
\begin{equation}\label{eq:23july16}
E((0,0,\hat{p}^a_3), (0,0,\hat{p}^b_3))=\min \{E((0,0,q^a_3), (0,0, q^b_3)): (q^a_3, q^b_3)\in P^*\},
\end{equation}
 and $\left(\psi^a_{(\hat{p}^a, \hat{p}^b)}, \psi^b_{(\hat{p}^a, \hat{p}^b)}\right)$ is the unique solution of (\ref{eq:19}). Moreover
 \begin{equation}\label{eq:24july16}
 \displaystyle\lim_n \frac{S_n((p^a_n, p^b_n))}{(h^a_n)^2}=E((0,0,\hat{p}^a_3), (0,0,\hat{p}^b_3)).
\end{equation}
\end{teo}

\vspace{0.5 cm}

\begin{teo}\label{teo:2july16}
Assume (\ref{eq:1}) with $\ell=0$ and (\ref{eq:14}). For every $n\in \NN$, let $(p^a_n, p^b_n)$ a solution of $(\ref{eq:13july16})$, and $\left(\phi^a_{(p^a_n, p^b_n)}, \phi^b_{(p^a_n, p^b_n)}\right)$ be the unique solution of $(\ref{eq:11})$ with $(p^a, p^b)=(p^a_n, p^b_n)$. Moreover, let $E_0$ defined by  (\ref{eq:18*}) and (\ref{eq:19*}). Then there exist an increasing sequence of positive integer numbers $\{n_i\}_{i\in\NN}$ and (in possible dependence on the subsequence) $\hat{p}^a\in H^1(]0,1[)$ such that (\ref{eq:20*})-(\ref{eq:22*}) hold, where $\hat{p}^a$ solves
\begin{equation}\label{eq:23*july16}
E_0(\hat{p}^a)=\min \{E_0(q^a): q^a\in H^1(]0,1[)\},
\end{equation}
 and $\psi^a_{\hat{p}^a}$ is the unique solution of (\ref{eq:19*}) with $q^a=\hat{p}^a$. Moreover
 \begin{equation}\label{eq:24*july16}
 \displaystyle\lim_n \frac{S_n((p^a_n, p^b_n))}{(h^a_n)^2}=E_0(\hat{p}^a).
\end{equation}
\end{teo}

\vspace{0.5 cm}

\begin{teo}\label{teo:3july16}
Assume (\ref{eq:1}) with $\ell=+\infty$ and $h^b_n\ll \sqrt{h^a_n}$, and (\ref{eq:14}). For every $n\in \NN$, let $(p^a_n, p^b_n)$ a solution of $(\ref{eq:13july16})$, and $\left(\phi^a_{(p^a_n, p^b_n)}, \phi^b_{(p^a_n, p^b_n)}\right)$ be the unique solution of $(\ref{eq:11})$ with $(p^a, p^b)=(p^a_n, p^b_n)$. Moreover, let  $E_{\infty}$ defined by (\ref{eq:18***}), (\ref{eq:19***}). Then there exist an increasing sequence of positive integer numbers $\{n_i\}_{i\in\NN}$ and (in possible dependence on the subsequence) $\hat{p}^b_3\in H^1(\Theta)$ such that
\begin{equation}\label{eq:20***july16}
\begin{array}{rcl}
          \displaystyle\frac{h^a_n}{\sqrt{h^b_n}}p^a_{n_i}\rightarrow 0 & \text{strongly in} & \left(H^1(\Omega^a)\right)^3\,  \text{ and strongly in}\,  \left(L^4(\Omega^a)\right)^3,\\
          p^b_{n_i}\rightarrow (0,0,\hat{p}^b_3) & \text{strongly in} & \left(H^1(\Omega^b)\right)^3\,  \text{ and strongly in}\,  \left(L^4(\Omega^b)\right)^3,
\end{array}
\end{equation}
and (\ref{eq:21***})-(\ref{eq:22***}) hold and
where $\hat{p}^b_3$ solves
\begin{equation}\label{eq:23***july16}
E_{\infty}(\hat{p}^b_3)=\min \{E_{\infty}(q^b_3):  q^b_3\in H^1(\Theta)\},
\end{equation}
 and $\psi^b_{\hat{p}^b}$ is the unique solution of (\ref{eq:19***}) with $q^b=\hat{p}^b$. Moreover
 \begin{equation}\label{eq:24***july16}
 \displaystyle\lim_n \frac{S_n((p^a_n, p^b_n))}{h^b_n}=E_{\infty}(\hat{p}^b_3).
\end{equation}
\end{teo}

\vspace{0.5 cm}

\begin{teo}\label{teonewjuly17}
Assume (\ref{eq:1}) with $\ell\in[0,+\infty]$ and (\ref{eq:14}). For every $n\in \NN$, let $(p^a_n, p^b_n)$ a solution of $(\ref{eq:13july16new})$, and $\left(\phi^a_{(p^a_n, p^b_n)}, \phi^b_{(p^a_n, p^b_n)}\right)$ be the unique solution of $(\ref{eq:11})$ with $(p^a, p^b)=(p^a_n, p^b_n)$. We have the convergences in Theorem \ref{teo:1}, Theorem \ref{teo:2} and Theorem \ref{teo:3} for the cases $q\in]0,+\infty[$, $q=0$ and $q=\infty$ respectively. Moreover,
 \begin{align}
 &\label{neweqjuly17911}
 \displaystyle\lim_n \frac{S_n((p^a_n, p^b_n))}{(h^a_n)^2}= \displaystyle\lim_n \frac{E_n((p^a_n, p^b_n))}{(h^a_n)^2}=E((0,0,\hat{p}^a), (\hat{p}^b_1, \hat{p}^b_2,0)),\, (\hat{p}^a, \hat{p}^b_1, \hat{p}^b_2)\in P,\,\text{if}\, q\in ]0,+\infty[,\\
&\label{neweqjuly17912}
 \displaystyle\lim_n \frac{S_n((p^a_n, p^b_n))}{(h^a_n)^2}= \displaystyle\lim_n \frac{E_n((p^a_n, p^b_n))}{(h^a_n)^2}=E_0(\hat{p}^a),\quad \hat{p}^a\in P_0 \quad \text{if}\,\, q=0,\\
&\label{neweqjuly17913}
 \displaystyle\lim_n \frac{S_n((p^a_n, p^b_n))}{h^b_n}= \displaystyle\lim_n \frac{E_n((p^a_n, p^b_n))}{h^b_n}=
 E_{\infty}(\hat{p}^b_3),\quad \hat{p}^b_3\in P_{\infty}\quad \text{if}\,\, q=+\infty.
\end{align}
\end{teo}

\section*{Acknowledgments}
This work was inspired by the great contributions of Prof. Antonio Gaudiello to the theory of thin structures. P. Hern\'andez-Llanos was supported by Agencia Nacional de Investigaci\'on y Desarro\-llo, FONDECYT Postdoctorado 2023 grant no. 3230202. The warm hospita\-lity at the Instituto de Ciencias de la Ingenier\'ia of the Universidad de O'Higgins is gratefully acknowledged by P. Hern\'andez-Llanos.


\vspace{1.5cm}

\end{document}